\RequirePackage{fix-cm}
\documentclass[smallextended]{svjour3}       % onecolumn (second format)
\smartqed  % flush right qed marks, e.g. at end of proof
\usepackage{graphicx}
\usepackage{amssymb}
\usepackage[colorlinks=true,linkcolor=blue]{hyperref}
\usepackage{amsmath,amssymb}
\usepackage{graphicx,graphics}% Include figure files
\usepackage{subfigure}
\usepackage{dcolumn}% Align table columns on decimal point
\usepackage{bm}% bold math
\usepackage{float}

\newcommand{\cQ}{\mathcal{Q}}

\renewcommand{\div}{\mbox{\rm div\,}}
\newcommand{\curl}{\mbox{\rm curl\,}}

\newcommand{\cM}{\mathcal{M}}
\newcommand{\cV}{\mathcal{V}}
\newcommand{\cW}{\mathcal{W}}

\newcommand{\R}{\mathbf{R}}

\usepackage[misc]{ifsym}
\usepackage{latexsym}
\begin{document}

\title{Unconditionally energy stable numerical schemes for the
three-dimensional magneto-micropolar equations%\thanks{Grants or other notes
%about the article that should go on the front page should be
%placed here. General acknowledgments should be placed at the end of the article.}
}
%\subtitle{Well-posedness and finite element approximation for MHD with temperature-dependent parameters }

\titlerunning{Unconditionally energy stable schemes for the
three-dimensional magneto-micropolar equations}        % if too long for running head

\author{Hailong Qiu \textsuperscript{1} %\and Liquan Mei
}

\authorrunning{H. Qiu} % if too long for running head

\institute{Hailong Qiu \at
School of Mathematics and Physics, Yancheng Institute of Technology,
Yancheng, 224051, China.\\
 This work is supported by the Natural Science Foundation of China (No. 11701498).\\
              %\email{hail$\underline{\phantom{a}}$qiu@126.com}
              \email{hail\_qiu@126.com}           %  \\
%             \emph{Present address:} of F. Author  %  if needed
}

\date{Received: date / Accepted: date}
% The correct dates will be entered by the editor

\maketitle

\begin{abstract}
In this paper we consider unconditionally energy stable numerical schemes for
the nonstationary $3$D magneto-micropolar equations
that describes the microstructure of rigid microelements in
electrically conducting fluid flow under some magnetic field.
The first scheme is comprised of the Euler semi-implicit discretization
in time and conforming finite element/stabilized finite element in space.
The second one is based on Crank-Nicolson discretization
in time and extrapolated treatment of the nonlinear terms such that skew-symmetry
properties are retained.
 We prove that the proposed schemes are
unconditionally energy stable.  Some error estimates for the velocity field, the magnetic field,
 the micro-rotation field and  the
fluid pressure are obtained. Furthermore, we establish some
 first-order decoupled numerical schemes.
 Numerical tests are provided to
check the theoretical  rates and unconditionally energy stable.
% \PACS{PACS code1 \and PACS code2 \and more}
\end{abstract}

\section{Introduction}
\label{intro}
Magneto-micropolar system  is used to describe the
motion of an incompressible, electrically conducting micropolar fluid under the action of magnetic fields.
The magneto-micropolar equations are governed by a combination of magneto-hydrodynamics equations and angular momentum equations.
The research of magneto-micropolar problems is a importance in
mathematical theory and practical applications, such as effective design, control and fabrication for micro-fluidic models \cite{Arkilic1997,Holmes1968,Wu1983,Wilding1994}.

  In this paper, we are interested in approximating the solution of the nonstationary magneto-micropolar fluid flows. The governing equations are defined by
\begin{subequations}\label{eq1.1}
\begin{alignat}{2} \label{eq1.1a} \partial_t \textbf{u}-(\nu+\nu_r)\Delta\textbf{u}
+(\textbf{u}\cdot \nabla)\textbf{u}
+S\textbf{B}\times&\mbox{curl}\textbf{B}+\nabla p
\\\qquad\nonumber&=\nu_r\mbox{curl}\textbf{w}+\textbf{f},  &&\ \mbox{in}\  \Omega_T,\\
\partial_t \textbf{B}+\mu\mbox{curl}(\mbox{curl} \textbf{B})-\mbox{curl}(\textbf{u}\times \textbf{B})&=\textbf{0}, &&\ \mbox{in}\  \Omega_T,\\
\mbox{div} \textbf{u}&=0, &&\ \mbox{in}\  \Omega_T,\\
\mbox{div} \textbf{B}&=0, &&\ \mbox{in} \ \Omega_T,\\
 \partial_t \textbf{w}-(c_a+c_d)\Delta\textbf{w}+(\textbf{u}\cdot\nabla)\textbf{w}+2\nu_r\textbf{w}&=\nu_r\mbox{curl}\textbf{u}\\\nonumber
 +(c_0+c_d-c_a)\nabla(\nabla\cdot\textbf{w})&+\textbf{g},  &&\ \mbox{in} \ \Omega_T,
\end{alignat}
\end{subequations}
complemented with the following homogeneous boundary conditions
\begin{subequations}\label{eq2.1}
\begin{alignat}{2} \label{eq2.1a}
 \textbf{u}&=\textbf{0}, \quad &&\ \mbox{on} \ \partial\Omega_T,\\
 \textbf{B}\cdot \textbf{n}&=0, \quad &&\ \mbox{on} \ \partial\Omega_T,\\
\mbox{curl \textbf{B}}\times \textbf{n}&=\textbf{0}, \quad &&\ \mbox{on} \  \partial\Omega_T,\\
\textbf{w}&=\textbf{0},\ \quad &&\ \mbox{on}\  \partial\Omega_T,
\end{alignat}
\end{subequations}
 and initial conditions
\begin{align}\label{eq3}
\textbf{u}(\textbf{x},0)=\textbf{u}_0(\textbf{x}),\ \
\textbf{B}(\textbf{x},0)=\textbf{B}_0(\textbf{x}), \ \
\textbf{w}(\textbf{x},0)=\textbf{w}_0(\textbf{x}),\ \ \mbox{in}\ \ \Omega,
\end{align}
 where $\Omega_T:=(0,T]\times\Omega$, $T$ stands for time, and $\Omega \subset \R^3$ is a
bounded domain with Lipschitz continuous boundary $\partial\Omega$.
$\textbf{u}$, $p$, $\textbf{B}$ and $\textbf{w}$ are
the velocity field, the pressure, the magnetic field and the micro-rotation field,  respectively.
% $\mathbb{D}(\textbf{u}):=\frac{1}{2}(\nabla \textbf{u}+\nabla \textbf{u}^T)$ is the strain-rate tensor.
 $\textbf{f}$ and $\textbf{g}$ denote the body force term and the moment term, respectively.
 The positive constants $\nu$, $\nu_r$ and $\mu$ are
respectively the fluid viscosity, the micro-rotation viscosity and  the magnetic viscosity.
 The positive constants $c_0$, $c_a$ and $c_d$ denotes viscosity coefficients which satisfy $c_0+c_d>c_a$.
$S$ represents the coupling coefficient.

If the magnetic field $\textbf{B}= \textbf{0}$, the magneto-micropolar equations \eqref{eq1.1}-\eqref{eq3} becomes the usual micropolar problem \cite{Eringen1966}.
The micropolar fluid flows describes the behavior of the local rotational motion of fluid molecules included in a given fluid volume element.
  The field of micropolar problems attracted the attention of many scholars to research both
mathematical theory and numerical calculation \cite{1977Galdi,Dong2010,Chen2012,chenm2012,Nochetto14,Salgado15,CrnjaricZica2016,Ma2018,Slayi2021}.
In \cite{1977Galdi}, Galdi and Rionero discussed the global existence and uniqueness of a weak solution for  micropolar equations.
In \cite{chenm2012},  Chen and Miao considered the global well-posedness of the strong solution with small initial data for  micropolar equations.
In \cite{Salgado15}, Salgado studied the convergence analysis of fractional time-stepping schemes for incompressible  micropolar problem.
In \cite{Slayi2021}, Slayi and cooperators proposed stabilized gauge uzawa scheme and established some error estimates for an evolutionary
micropolar fluid flows.

 Magneto-micropolar fluid flows was firstly studied by Ahmadi and Shahinpoor in \cite{Ahmadi1974}, which was established by the classical magneto-hydrodynamics augmented with
 additional equations that explain conservation of angular momentum due to consideration of the micro-structure
in the fluid \cite{Chen2011,Petrosyan}.
As the mathematical importance and physical applications of the magneto-micropolar fluid flows, this
 system has been studied by many researchers \cite{1997Rojas,1998Rojas,Yamazaki2015,2016Silva,MaS2018,2019tan,2020Ravindran,2021Wang,Ai2022,2022Fan,2022Qin,Zou2022}. % More attention is focused on the
%mathematical theory and numerical computation of the magneto-micropolar equations.
 In \cite{1997Rojas}, Rojas-Medar applied spectral Galerkin method to establish the well-posedness of strong solutions for the magneto-micropolar fluid flows.
  In \cite{1998Rojas}, Rojas-Medar and Boldrini considered
the existence of weak solutions  for the nonstationary magneto-micropolar problem.
Recently, Ravindran in \cite{2020Ravindran} discussed unconditionally long time stable and convergence
of a BDF2 time-stepping scheme for nonstationary magneto-micropolar fluid
flows.

Developing efficient, accurate, and robust numerical schemes for solving the
multiple physical coupled system involved with incompressible fluid flows remains a great challenge
in $3$D magneto-micropolar model. There are many time-stepping methods of multiple physical
coupled problems at present \cite{2015Chen,2019Chen,2022Chen,Diegel,He2015,2017Hasler,L2013,2016Liu,Prohl2008,Tone2009,2022Wang}.
In \cite{L2013},  Layton studied two partitioned time-stepping algorithms for a reduced magneto-hydrodynamics
problem. In \cite{He2015}, He considered  an Euler semi-implicit scheme and obtained some
convergence for the magneto-hydrodynamics equations.
In \cite{2017Hasler}, Heister investigated a decoupled higher order
scheme for magneto-hydrodynamics equations.
Recently, in \cite{2022Wang}, the authors studied a Crank-Nicolson time-stepping projection
scheme and derived some optimal error estimates for magnetohydrodynamic problem.

In this paper, we propose some unconditionally energy stable
 numerical schemes for the  nonstationary $3$D magneto-micropolar equations.
The first scheme is comprised of the Euler semi-implicit discretization
in time and conforming finite element/stabilized finite element in space.
The second one is based on Crank-Nicolson discretization
in time and extrapolated treatment of the nonlinear terms such that skew-symmetry
properties are  preserved.
 We prove that the proposed schemes are
unconditionally energy stable. Some convergence rates of the velocity field, the magnetic field,
the micro-rotation  field and the
fluid pressure are established. In addition, we propose some
fully discrete first-order decoupled numerical schemes.
Numerical tests are given to
conform the theoretical results of our proposed schemes.

The remainder of this paper is organized as follows.
In section 2, we introduce some preliminaries
and the weak formulation
 of magneto-micropolar equations.  In section 3, we propose
a first-order numerical scheme and obtain its unconditionally stable.
In section 4, we derive some error estimates for the velocity field, the magnetic field,
the micro-rotation field and the pressure. In section 5, we set up
a second-order numerical scheme and give some unconditionally stable and error estimates.
Moreover, we establish some first-order decoupled  numerical schemes.
 In section 6, we give
 several numerical rests to check the theoretical results.

\section{\label{Sec2} Continuous problem}

For $k\geq 0$ and $p\in [1,+\infty]$, we denote by $L^p(\Omega)$ and $W^{k,p}(\Omega)$
the usual Lebesgue spaces and standard Sobolev spaces equipped
with the norms  $\|\cdot\|_{L^p}$ and $\|\cdot\|_{W^{k,p}}$, respectively. In particular, let us denote $\|\cdot\|_{L^2}:=\|\cdot\|$.  We write $H^k(\Omega)$ in place of
the Sobolev space $W^{k,2}(\Omega)$ when $p = 2$,
  and denote the corresponding norm by $\|\cdot\|_{k}$.
%We denote the time discrete space $l^p(Z)$ associated with $L^p(0,T; Z)$.

To give a weak formulation of problem \eqref{eq1.1}-\eqref{eq3}, we introduce some special spaces
\begin{align*}
&\cV:={H}_0^1(\Omega)^3=\bigl\{\textbf{v}\in H^1(\Omega)^3: \textbf{v}|_{\partial\Omega}=0\bigr\}, \\
& \cW:={H}^1_n(\Omega)^3:=\bigl\{\textbf{w}\in H^1(\Omega)^3: \,\textbf{w}\cdot\textbf{n}|_{\partial\Omega}=0 \bigr\}, \\
& \cW_0:=\bigl\{\textbf{w}\in H^1(\Omega)^3: \,\curl\textbf{w}\times\textbf{n}|_{\partial\Omega}=0 \bigr\}, \\
&\cV_0:=\bigl\{\textbf{v}\in  \cV: \mbox{div}\textbf{v}=0\bigr\}, \\
& \cW_n: =\bigl\{\textbf{w}\in  \cW: \mbox{div}\textbf{w}=0 \bigr\},\\
%&\cQ:=\bigl\{\varphi\in H^1(\Omega)^3: \varphi|_{\partial\Omega}=0\bigr\}, \\
& \cM: =L_0^2(\Omega)=\bigl\{ q\in L^2(\Omega), \int_\Omega q d \textbf{x}=0 \bigr\}.
\end{align*}

%It is well known that
In the sequel, we will make use of the following inequalities  \cite{Temam1983,Girault1986,1991Gunzburger}.
\begin{subequations}\label{eq4.1}
\begin{alignat}{2} \label{eq4.1a}
\|\textbf{v}\|&\leq C\|\nabla\textbf{v}\|, \ \  &&\ \forall \ \textbf{v} \in \cV,\\\label{eq4.1b}
\|\textbf{v}\|_{L^6}&\leq   C\|\nabla\textbf{v}\|, \ \ &&\ \forall \ \textbf{v} \in \cV,\\\label{eq4.1c}
%\|\mathbb{D}\textbf{v}\|&\geq  c^{\star}\|\nabla\textbf{v}\|, \ \ &&\ \forall \ \textbf{v} \in \cV,\\\label{eq4.1d}
\|\textbf{v}\|_{L^\infty}&\leq   C\|\textbf{v}\|_1^{1/2}\|\textbf{v}\|^{1/2}_{2},  &&\ \forall\ \textbf{v} \in {H}^2(\Omega)^3,\\\label{eq4.1e}
\|\textbf{v}\|_{L^3}&\leq  C\|\textbf{v}\|^{1/2}\|\textbf{v}\|^{1/2}_{1},\ &&\ \forall\ \textbf{v}\in {H}^1(\Omega)^3,\\\label{eq4.1f}
\|\curl \textbf{B}\|^2+\|\div \textbf{B}\|^2&\geq c_{\star}  \|\nabla\textbf{B}\|^2,&&\ \forall\ \textbf{B} \in \cW,\\\label{eq4.1g}
\|\curl \textbf{B}\|\leq\sqrt{2}\|\nabla\textbf{B}\|,  \|\div\textbf{B}\|&\leq \sqrt{d}\|\nabla\textbf{B}\|, &&\ \forall\ \textbf{B} \in {H}^1(\Omega)^3,
\end{alignat}
\end{subequations}
where $C$ and $c_{\star}$ are positive constants depend on $\Omega$.

For simplicity, in what follows we define the bilinear terms
\begin{align*}
a_f(\textbf{u},\textbf{v})&=(\nu+\nu_r)\int_\Omega\nabla\textbf{u}\cdot\nabla\textbf{v}d\textbf{x},\\
 a_w(\textbf{w},\varphi)&=(c_a+c_d)\int_\Omega\nabla\textbf{w}\cdot\nabla\varphi d\textbf{x}+(c_0+c_d-c_a)\int_\Omega\nabla\cdot\textbf{w}\nabla\cdot\varphi d\textbf{x}, \\
a_{B}(\textbf{B},\textbf{H})&=\mu\int_\Omega\curl\textbf{B}\cdot\curl\textbf{H}d\textbf{x}
+\mu\int_\Omega\div\textbf{B}\cdot\div\textbf{H}d\textbf{x},\\
 d(\textbf{v},q)&=\int_\Omega q\mbox{div}\textbf{v}d\textbf{x},\qquad e(\textbf{u},\textbf{v})=\nu_r\int_\Omega\curl\textbf{u}\cdot\textbf{v}d\textbf{x},
\end{align*}
and trilinear terms
\begin{align*}
b(\textbf{w},\textbf{u},\textbf{v})&= \frac{1}{2}\int_\Omega[(\textbf{w}\cdot\nabla)\textbf{u}]\cdot \textbf{v}
 -[(\textbf{w}\cdot\nabla)\textbf{v}]\cdot \textbf{u} d\textbf{x}=\int_\Omega[(\textbf{w}\cdot\nabla)\textbf{u}]\cdot \textbf{v}
 +\frac{1}{2}[(\nabla\cdot\textbf{w})\textbf{u}]\cdot \textbf{v} d\textbf{x},\\
c_w(\textbf{w},\psi,\varphi)&= \frac{1}{2}\int_\Omega(\textbf{w}\cdot\nabla) \psi \varphi
 -(\textbf{w}\cdot\nabla)\varphi \psi  d\textbf{x}=\int_\Omega(\textbf{w}\cdot\nabla)\psi \varphi
 +\frac{1}{2}(\nabla\cdot\textbf{w})\varphi  \psi d\textbf{x},\\
c_{\widehat{B}}(\textbf{H},\textbf{B},\textbf{v})&=\int_\Omega \textbf{H}\times \mbox{curl}\textbf{B}\cdot\textbf{v}d\textbf{x},\
\ \quad c_{\widetilde{B}}(\textbf{u},\textbf{B},\textbf{H}) =\int_\Omega(\textbf{u}\times \textbf{B})\cdot\curl\textbf{H}d\textbf{x}.
\end{align*}
Additionally, using the definitions of $b(\cdot,\cdot,\cdot)$ and $c_w(\cdot,\cdot,\cdot)$, and thanks to integration by parts, we get
\begin{align}\label{eq5}
b(\textbf{u},\textbf{v},\textbf{v})&=0,\  &&\textbf{u}\in \cV, \textbf{v}\in {H}^1(\Omega)^3,\\\label{eq6}
c_w(\textbf{u},\psi,\psi)&=0,&&\textbf{u}\in \cV, \psi\in H^1(\Omega)^3.
\end{align}
Due to $(\textbf{B}\times \mbox{curl}\textbf{H},\textbf{v})=(\textbf{v}\times \textbf{H}, \mbox{curl}\textbf{B})$,  we arrive at
\begin{align}\label{eq7}
c_{\widehat{B}}(\textbf{B},\textbf{B},\textbf{u})-c_{\widetilde{B}}(\textbf{u},\textbf{B},\textbf{B})=0,\ \textbf{u}\in \cV,  \textbf{B} \in \cW.
\end{align}
In addition, the bilinear term  $d(\cdot,\cdot)$  satisfies the inf-sup condition \cite{Girault1986,Temam1983} as follows
\begin{align}\label{eq8}
\sup_{\textbf{v}\in \cV,\textbf{v}\neq \textbf{0}}\dfrac{d(\textbf{v},q)}{\|\textbf{v}\|_{1}}\geq \beta\|q\|,\ \ \  \forall q \in \cM,
\end{align}
where $\beta$ is a positive constant depending on $\Omega$.

The weak formulation for problem \eqref{eq1.1}-\eqref{eq3} can be read as: find $(\textbf{u}, p,\textbf{B},\textbf{w})\in \cV\times \cM \times \cW\times \cV$
 such that
 \begin{subequations}\label{eq9.1}
\begin{alignat}{2} \label{eq9.1a} (\partial_t\textbf{u},\textbf{v})+a_f(\textbf{u},\textbf{v})
+b(\textbf{u},\textbf{u},\textbf{v})+Sc_{\widehat{B}}(\textbf{B},\textbf{B},\textbf{v})
\\\qquad\qquad\nonumber-d(\textbf{v},p)&=e(\textbf{w},\textbf{v})+(\textbf{f},\textbf{v}), \\
d(\textbf{u},q)&=0, \\
(\partial_t\textbf{B},\textbf{H})+a_B(\textbf{B},\textbf{H})-c_{\widetilde{B}}(\textbf{u},\textbf{B},\textbf{H})&=0, \\
(\partial_t\textbf{w},\psi)+a_w(\textbf{w},\psi)+c_w(\textbf{u},\textbf{w},\psi)+2\nu_r(\textbf{w},\psi)&=e(\textbf{u},\psi)+(\textbf{g},\psi),
\end{alignat}
\end{subequations}
for all $(\textbf{v},q,\textbf{H},\psi)\in \cV\times \cM \times \cW\times \cV$.

 The results in respect to the well-posedness of weak solutions to the problem
 \eqref{eq1.1}-\eqref{eq3} in a slightly different setting were established in \cite{1998Rojas}.

\begin{theorem} The weak formulation \eqref{eq9.1} has at least one solution $(\textbf{u}, p,  \textbf{B}, \textbf{w})$  such that
  \begin{align}\label{eq10}
 \textbf{u}&\in L^{\infty}(0,T,L^2(\Omega)^d)\cap L^2(0,T,\cV), \ p\in L^2(0,T,\cM), \\\nonumber
 \textbf{B}&\in L^{\infty}(0,T,L^2(\Omega)^d)\cap L^2(0,T,\cW), \
\textbf{w}\in L^{\infty}(0,T,L^2(\Omega)^d)\cap L^2(0,T,\cQ).
 \end{align}
Moreover, problem \eqref{eq9.1}  has an unique weak solution in two spatial dimension.% ($d = 2$).
\end{theorem}

\section{\label{Sec3} First-order numerical schemes }

In this section, using the backward Euler semi-implicit discretization
in time and conforming finite element/stabilized finite element in space,
we establish a fully discrete first-order numerical scheme
for problem \eqref{eq9.1} and obtain some unconditionally energy stable results.

Suppose that $\mathbb{K}_h$ is a quasi-uniform
triangulation of the polygonal/polyhedral bounded domain $\Omega\subset \R^3$ into tetrahedra.
Moreover, we consider two kinds of finite element spaces:

\textbf{Case I}
 \begin{align*}
\cV^1_h&=\bigl\{\textbf{v}\in \textbf{C}^0(\overline{\Omega})^3\cap \cV:\textbf{v}|_{\mathbb{K}}\in {P}_{1b}(\mathbb{K})^3,\forall\ \mathbb{K}\in \mathbb{K}_h\bigr\}, \\
 \cM_h&=\bigl\{q\in C^0(\overline{\Omega})\cap\cM:q|_{\mathbb{K}}\in P_1(\mathbb{K}),\forall\ \mathbb{K}\in \mathbb{K}_h\bigr\}, \\
\cW^1_h&=\bigl\{\textbf{w}\in \textbf{C}^0(\overline{\Omega})^3\cap \cW:\textbf{w}|_{\mathbb{K}}\in {P}_{1b}(\mathbb{K})^3,\forall\ \mathbb{K}\in \mathbb{K}_h\bigr\}, \\
 \cQ^1_h&=\bigl\{\omega\in \textbf{C}^0(\overline{\Omega})^3\cap\cQ:\omega|_{\mathbb{K}}\in P_{1b}(\mathbb{K})^3,\forall\ \mathbb{K}\in \mathbb{K}_h\bigr\},
 \end{align*}
where
$${P}_{1b}(\mathbb{K})=\bigl\{{v}_h\in {C}^0(\overline{\Omega}): {v}_h|_\mathbb{K} \in {P}_1(\mathbb{K})\oplus \mbox{span}\{\hat{b}\}, \ \mathbb{K}\in \mathbb{K}_h\bigr\},$$
$\hat{b}$ stands for a bubble function, and $P_1(\mathbb{K})$ presents the set
of polynomials with degree less than equal to $1$ over $\mathbb{K}$.

\textbf{Case II}
 \begin{align*}
\cV^2_h&=\bigl\{\textbf{v}\in \textbf{C}^0(\overline{\Omega})^3\cap \cV:\textbf{v}|_{\mathbb{K}}\in {P}_{1}(\mathbb{K})^3,\forall\ \mathbb{K}\in \mathbb{K}_h\bigr\}, \\
 \cM_h&=\bigl\{q\in C^0(\overline{\Omega})\cap\cM:q|_{\mathbb{K}}\in P_1(\mathbb{K}),\forall\ \mathbb{K}\in \mathbb{K}_h\bigr\}, \\
\cW^2_h&=\bigl\{\textbf{w}\in \textbf{C}^0(\overline{\Omega})^3\cap \cW:\textbf{w}|_{\mathbb{K}}\in {P}_{1}(\mathbb{K})^3,\forall\ \mathbb{K}\in \mathbb{K}_h\bigr\}, \\
 \cQ^2_h&=\bigl\{\omega\in \textbf{C}^0(\overline{\Omega})^3\cap\cQ:\omega|_{\mathbb{K}}\in P_{1}(\mathbb{K})^3,\forall\ \mathbb{K}\in \mathbb{K}_h\bigr\}.
 \end{align*}

%For case I, we assume that the discrete inf-sup condition \cite{Temam1983,Girault1986} holds for stable mixed finite element pair $(\cV^1_h,\cM_h)$.
For case I, the mixed finite element pair $(\cV^1_h,\cM_h)$ is stable and satisfy the following inf-sup condition \cite{Temam1983,Girault1986}.

\textbf{Assumption A1:} The following inequality holds
\begin{align*}
\sup_{\textbf{v}_h\in\cV^1_h,\textbf{v}_h\neq 0}\dfrac{d(\textbf{v}_h,q_h)}{\|\nabla\textbf{v}_h\|}\geq {\beta}_0\|q_h\|, \ \forall q_h \in {\cM}_h,
 \end{align*}
where ${\beta}_0$ is a positive constant independent of $h$.

For case II, as is noted that the finite element pair $(\cV^2_h,\cM_h)$ is not stable. Thus, let us denote
 \begin{equation*}\overline{\mathbb{B}}(\textbf{u}_h,p_h;\textbf{v}_h,q_h):=a_f(\textbf{u}_h,v_h)-d(\textbf{v}_h,p_h)+d(\textbf{u}_h,q_h)+
G(p_h,q_h),\end{equation*}
where the stabilized term $G(p_h,q_h)$ is defined by
$$G(p_h,q_h):= ((I-\Pi_1) p_h,(I-\Pi_1) q_h),$$
here $\Pi_1$ is a $L^2$-projection operator.

%As is known that the $(\cV_h,\cM_h)$ is a stable mixed finite element pair.
\begin{lemma}\cite{Qiu2018} For all $(\textbf{u}_h,p_h),(\textbf{v}_h, q_h) \in \cV^2_h\times \cM_h$, $\overline{\mathbb{B}}(\cdot,\cdot;\cdot,\cdot)$ satisfies the continuity property
\begin{equation*}| \bar{\mathbb{B}}(\textbf{u}_h,p_h;\textbf{v}_h,q_h)|\leq
\beta_1\bigl(\|\nabla \textbf{u}\|+\| p_h\|\bigr)\bigl(\|\nabla\textbf{v}_h\|+\| q_h\|\bigr),\quad \forall \
(\textbf{u}_h,p_h),(\textbf{v}_h,q_h)\in (\cV^2_h,\cM_h),\end{equation*}
and the weak coercivity property
\begin{equation*}\beta_2(\|\nabla\textbf{u}_h\|+\| p_h\|)\leq
\sup_{(v,q)\in(\widetilde{\cV^2_h},\cM_h)} \frac{
\bar{\mathbb{B}}(\textbf{u}_h,p_h;\textbf{v}_h,q_h)}{\|\nabla \textbf{v}_h\|+\| q_h\|},\quad \forall \
(\textbf{u}_h,p_h)\in (\widetilde{\cV^2_h},\cM_h),\end{equation*} where $ \widetilde{\cV^2_h}
=\bigl\{v\in \cV_0: \ v|_{\mathbb{K}}\in P_1(\mathbb{K})^d, \ \forall \ \mathbb{K}\in \mathbb{K}_h\bigr\}$ is the
finite element subspace of $\cV_0$, $\beta_1$ and $\beta_2$ are two positive
constants independent of $h$.
\end{lemma}

 Let $N$ be a positive integer, $\Delta t :=\frac{T}{N}$ and $t_n=n\Delta t$ for $n=0,1,\cdots,N$. %{\color{blue}
  Let ${\phi}^n={\phi}(t_n)$ and
denote $\mathcal{D}{\phi}^n:=\frac{\phi^n-\phi^{n-1}}{\Delta t}$.%}

 We now propose a fully discrete first-order  numerical scheme of problem \eqref{eq9.1} as follows.

\textbf{Scheme 3.1:}  Given $(\textbf{u}^0_h,\textbf{B}^0_h,\textbf{w}^0_h) \in \cV^1_h\times \cW^1_h\times \cQ^1_h$,
find $(\textbf{u}^n_h, p^n_h, \textbf{B}^n_h, \textbf{w}^n_h) \in {\cV}^1_h\times \cM_h\times \cW^1_h\times \cQ^1_h$, such that
 \begin{subequations}\label{eq11.1}
\begin{alignat}{2} \label{eq11.1a}
 (\mathcal{D}\textbf{u}^n_h,\textbf{v}_h)+a_f(\textbf{u}^n_h,\textbf{v}_h)
+b(\textbf{u}^{n-1}_h,\textbf{u}^n_h,\textbf{v}_h)\\\nonumber\quad
+Sc_{\widehat{B}}(\textbf{B}^{n-1}_h,\textbf{B}^n_h,\textbf{v}_h)-d(\textbf{v}_h,p^n_h)
&=e(\textbf{w}_h^n,\textbf{v}_h)+(\textbf{f}^n,\textbf{v}_h), \\
d(\textbf{u}^n_h,q_h)&=0,\\
(\mathcal{D}\textbf{B}^n_h,\textbf{H}_h)+a_B(\textbf{B}^{n}_h,\textbf{H}_h)
-c_{\widetilde{B}}(\textbf{u}^{n}_h,\textbf{B}^{n-1}_h,\textbf{H}_h)&=0, \\
(\mathcal{D}\textbf{w}^n_h,\psi_h)+a_w(\textbf{w}^{n}_h,\psi_h)+c_w(\textbf{u}^{n-1}_h,\textbf{w}^{n}_h,\psi_h)
\\\nonumber\quad+2\nu_r(\textbf{w}_h^n,\psi_h)&=e(\textbf{u}_h^n,\psi_h)+(\textbf{g}^n,\psi_h),
\end{alignat}
\end{subequations}
for all $(\textbf{v}_h,q_h,\textbf{H}_h,\psi_h)\in {\cV}^1_h\times \cM_h\times \cW^1_h\times \cQ^1_h$.

 Considering stabilized finite element pair, we present a new scheme as follows:

\textbf{Scheme 3.2:}  Given $(\textbf{u}^0_h,\textbf{B}^0_h,\textbf{w}^0_h) \in \cV^2_h\times \cW^2_h\times \cQ^2_h$,
find $(\textbf{u}^n_h, p^n_h, \textbf{B}^n_h, \textbf{w}^n_h) \in {\cV}^2_h\times \cM_h\times \cW^2_h\times \cQ^2_h$, such that
 \begin{subequations}\label{eq12.1}
\begin{alignat}{2} \label{eq12.1a}
 (\mathcal{D}\textbf{u}^n_h,\textbf{v}_h)+\overline{\mathbb{B}}(\textbf{u}^n_h,p^n_h;\textbf{v}_h,q_h)
+b(\textbf{u}^{n-1}_h,\textbf{u}^n_h,\textbf{v}_h)\\\nonumber\quad
+Sc_{\widehat{B}}(\textbf{B}^{n-1}_h,\textbf{B}^n_h,\textbf{v}_h)
&=e(\textbf{w}_h^n,\textbf{v}_h)+(\textbf{f}^n,\textbf{v}_h), \\
(\mathcal{D}\textbf{B}^n_h,\textbf{H}_h)+a_B(\textbf{B}^{n}_h,\textbf{H}_h)
-c_{\widetilde{B}}(\textbf{u}^{n}_h,\textbf{B}^{n-1}_h,\textbf{H}_h)&=0, \\
(\mathcal{D}\textbf{w}^n_h,\psi_h)+a_w(\textbf{w}^{n}_h,\psi_h)+c_w(\textbf{u}^{n-1}_h,\textbf{w}^{n}_h,\psi_h)
\\\nonumber\quad+2\nu_r(\textbf{w}_h^n,\psi_h)&=e(\textbf{u}_h^n,\psi_h)+(\textbf{g}^n,\psi_h),
\end{alignat}
\end{subequations}
for all $(\textbf{v}_h,q_h,\textbf{H}_h,\psi_h)\in  {\cV}^2_h\times \cM_h\times \cW^2_h\times \cQ^2_h$.

The following technique lemma is very important to analyse error estimate. % for $\partial_t\textbf{u}$, $\partial_t\textbf{B}$ and $\partial_t\textbf{w}$.
\begin{lemma}{\color{blue}\cite{Tabata2005,Qiu2020}}
For any $h>0$ and $\Delta t>0$, the following inequalities hold
  \begin{align*}
(i), b(\textbf{u}_h,\textbf{v}_h,\textbf{w}_h)&\leq C \min\bigl\{h^{-1/2},\Delta t^{-1/2}\bigr\}\|\nabla\textbf{u}_h\|\|\nabla\textbf{v}_h\|
  \\\nonumber&\quad\times\bigl(\|\textbf{w}_h\|+\sqrt{\Delta t}\|\nabla\textbf{w}_h\|\bigr), \forall \ \textbf{u}_h,\textbf{v}_h,\textbf{w}_h \in \cV_h, \\
(ii),  c_w(\textbf{u}_h,\psi_h,\phi_h)&\leq C \min\bigl\{h^{-1/2},\Delta t^{-1/2}\bigr\}\|\nabla\textbf{u}_h\|\|\nabla\psi_h\|
  \\\nonumber&\quad\times\bigl(\|\psi_h\|+\sqrt{\Delta t}\|\nabla\psi_h\|\bigr), \forall \ \textbf{u}_h \in \cV_h, \psi_h,\phi_h \in \cQ_h,\\
(iii), c_{\widehat{\textbf{B}}}(\textbf{H}_h,\textbf{B}_h,\textbf{v}_h)&\leq C \min\bigl\{h^{-1/2},\Delta t^{-1/2}\bigr\}\|\nabla\textbf{H}_h\|\|\nabla\textbf{B}_h\|
  \\\nonumber&\quad\times\bigl(\|\textbf{v}_h\|+\sqrt{\Delta t}\|\textbf{v}_h\|\bigr), \forall \ \textbf{H}_h,\textbf{B}_h \in \cW_h,\textbf{v}_h \in \cV_h,\\
(iv), c_{\widetilde{\textbf{B}}}(\textbf{u}_h,\textbf{B}_h,\textbf{H}_h)&\leq C \min\bigl\{h^{-1/2},\Delta t^{-1/2}\bigr\}\|\nabla\textbf{u}_h\|\|\nabla\textbf{B}_h\|
  \\\nonumber&\quad\times\bigl(\|\textbf{H}_h\|+\sqrt{\Delta t}\|\textbf{H}_h\|\bigr), \forall\  \textbf{u}_h \in \cV_h,\textbf{B}_h,\textbf{H}_h \in \cW_h.
 \end{align*}
\end{lemma}

We may now state and derive some unconditionally energy stable results for the fully discrete first-order numerical scheme  \eqref{eq11.1}.
\begin{theorem}  %Assume \textbf{Hypotheses A1-A3} hold.
Let $(\textbf{f},\textbf{g})\in L^2(0,T; H^{-1}(\Omega)^3)\times L^2(0,T; H^{-1}(\Omega)^3)$.
Then the finite element solution $(\textbf{u}^n_h,\textbf{B}^n_h,\textbf{w}^n_h)$ of problem \eqref{eq11.1} satisfy the following bounds
   \begin{align} \label{eq13} &\|\textbf{u}^N_h\|^2+S\|\textbf{B}^N_h\|^2
+\|\textbf{w}^N_h\|^2
\\\nonumber&\quad+{(\nu+\nu_r)}\Delta t\sum^{N}_{i=1}\|\nabla{\textbf{u}}^n_h\|^2+c_{\star}\mu\Delta tS\sum^{N}_{i=1}\|\nabla{\textbf{B}}^n_h\|^2\\\nonumber&\quad
+{\Delta t(c_a+c_d)}\sum^{N}_{i=1}\|\nabla{\textbf{w}}^n_h\|^2+(c_0+c_d-c_a)\Delta t\sum^{N}_{i=1}\|\nabla\cdot{\textbf{w}}^n_h\|^2\\\nonumber&
\leq C.
%\frac{\Delta t}{(\nu+\nu_r)}\sum^{N}_{i=1}\|\textbf{f}^n\|^2_{-1}+ \frac{\Delta t}{(c_a+c_d)}\sum^{N}_{i=1}\|\textbf{g}^n\|^2_{-1}+2\|\textbf{u}^{0}_h\|^2\\\nonumber&\quad
%+2S\|\textbf{B}^{0}_h\|^2
%+2\|\textbf{w}^{0}_h\|^2.
 \end{align}
\end{theorem}
\noindent \textit{Proof:}\quad Setting $(\textbf{v}_h, q_h, \textbf{H}_h,\psi_h)
=({\textbf{u}}^n_h,0,{\textbf{B}}^n_h,\textbf{w}^n_h)$ in \eqref{eq11.1}, using \eqref{eq5}-\eqref{eq7}, we get
\begin{align} \label{eq14}(\mathcal{D}{\textbf{u}}^n_h,{\textbf{u}}^n_h)&+S(\mathcal{D}{\textbf{B}}^n_h,{\textbf{B}}^n_h)
+(\mathcal{D}{\textbf{w}}^n_h,{\textbf{w}}^n_h)\\\nonumber&
+(\nu+\nu_r)\|\nabla{\textbf{u}}^n_h\|^2+c_{\star}\mu S\|\nabla{\textbf{B}}^n_h\|^2+2\nu_r\|{\textbf{w}}^n_h\|^2\\\nonumber&
+(c_a+c_d)\|\nabla{\textbf{w}}^n_h\|^2+(c_0+c_d-c_a)\|\nabla\cdot{\textbf{w}}^n_h\|^2\\\nonumber&
\leq e({\textbf{w}}^n_h,{\textbf{u}}^n_h)+e({\textbf{u}}^n_h,{\textbf{w}}^n_h)+(\textbf{f}^n,{\textbf{u}}^n_h)+(\textbf{g}^n,{\textbf{w}}^n_h).
 \end{align}
Applying H\"{o}lder inequality, Young's inequality and Cauchy inequality,  one finds that
\begin{align*}
e({\textbf{w}}^n_h,{\textbf{u}}^n_h)&\leq \sqrt{2}\nu_r\|\nabla{\textbf{u}}^n_h\|\|{\textbf{w}}^n_h\|\leq \frac{\nu_r}{8}\|\nabla{\textbf{u}}^n_h\|^2+4\nu_r\|{\textbf{w}}^n_h\|^2,\\
e({\textbf{u}}^n_h,{\textbf{w}}^n_h)&\leq \sqrt{2}\nu_r\|{\textbf{w}}^n_h\|\|\nabla{\textbf{u}}^n_h\|\leq \frac{\nu_r}{8}\|\nabla{\textbf{u}}^n_h\|^2+4\nu_r\|{\textbf{w}}^n_h\|^2,\\
(\textbf{f}^n,{\textbf{u}}^n_h)&\leq\frac{(\nu+\nu_r)}{4}\|\nabla{\textbf{u}}^n_h\|^2+\frac{1}{(\nu+\nu_r)}\|\textbf{f}^n\|^2_{-1},\\
(\textbf{g}^n,{\textbf{w}}^n_h)&\leq \frac{c_a+c_d}{2}\|\nabla{\textbf{w}}^n_h\|^2+\frac{1}{2(c_a+c_d)}\|\textbf{g}^n\|^2_{-1}.
 \end{align*}
 Inserting above inequalities into  \eqref{eq14} and using \eqref{eq4.1g}, we have
\begin{align} \label{eq15} (\mathcal{D}{\textbf{u}}^n_h,{\textbf{u}}^n_h)&+S(\mathcal{D}{\textbf{B}}^n_h,{\textbf{B}}^n_h)
+(\mathcal{D}\textbf{w}^n_h,\textbf{w}^n_h)\\\nonumber&
+\frac{(\nu+\nu_r)}{2}\|\nabla{\textbf{u}}^n_h\|^2+c_{\star}\mu S\|\nabla{\textbf{B}}^n_h\|^2\\\nonumber&
+\frac{(c_a+c_d)}{2}\|\nabla{\textbf{w}}^n_h\|^2+(c_0+c_d-c_a)\|\nabla\cdot{\textbf{w}}^n_h\|^2\\\nonumber&
\leq \frac{1}{2(\nu+\nu_r)}\|\textbf{f}^n\|^2_{-1} +\frac{1}{2(c_a+c_d)}\|\textbf{g}^n\|^2_{-1}+6\nu_r\|{\textbf{w}}^n_h\|^2.
 \end{align}
 Utilizing the elementary identity
 \begin{align*}
 2a(a-b)=a^2-b^2+(a-b)^2,\ \ \ \mbox{for all} \ a,\ b\in \R^3,
 \end{align*}
 we derive
 \begin{align} \label{eq16} &\|\textbf{u}^n_h\|^2+S\|\textbf{B}^n_h\|^2
+\|\textbf{w}^n_h\|^2+\|\textbf{u}^n_h-\textbf{u}^{n-1}_h\|^2+S\|\textbf{B}^n_h-\textbf{B}^{n-1}_h\|^2
\\\nonumber&\quad+\|\textbf{w}^n_h-\textbf{w}^{n-1}_h\|^2
+\frac{(\nu+\nu_r)}{2}\Delta t\|\nabla{\textbf{u}}^n_h\|^2+c_{\star}\mu\Delta tS\|\nabla{\textbf{B}}^n_h\|^2\\\nonumber&\quad
+\frac{(c_a+c_d)}{2}\|\nabla{\textbf{w}}^n_h\|^2+(c_0+c_d-c_a)\|\nabla\cdot{\textbf{w}}^n_h\|^2\\\nonumber&
\leq \frac{\Delta t}{2(\nu+\nu_r)}\|\textbf{f}^n\|^2_{-1}+ \frac{\Delta t}{2(c_a+c_d)}\|\textbf{g}^n\|^2_{-1}+\|\textbf{u}^{n-1}_h\|^2\\\nonumber&\quad
+S\|\textbf{B}^{n-1}_h\|^2
+\|\textbf{w}^{n-1}_h\|^2+6\nu_r\Delta t\|{\textbf{w}}^n_h\|^2.
 \end{align}
  Summing \eqref{eq16} from $n = 1$ to $N$ and using the discrete Gr\"{o}nwall Lemma, the desired result \eqref{eq13} holds.
 The proof is completed.
 $$\eqno\Box$$

  \textbf{Remark 1:} Applying similar lines as the proof of {Theorem 3.1},
 we can obtain a similar stability result for the scheme $3.2$, here we omit it.
%%%%%%%%%%%%%%%%%%%%%%%%%%%%%%%%%%%%%%%%%%%%%%%%%%%%%%%%%%%%%%%%%%%%%%%%%%%%%%%%%%%%%%%%%%%%%%%%%%%%%%%%%%%%%%%%%%%%

\section{\label{Sec4}Error analysis}

In this section, we derive the error estimates of the first-order numerical scheme \eqref{eq11.1}.
 We assume that the solution $(\textbf{u},p,\textbf{B},\textbf{w})$ of problem \eqref{eq9.1} has the following regularity result.

 \textbf{Assumpation A2:} The solution $(\textbf{u},p,\textbf{B},\textbf{w})$ of weak formulation \eqref{eq9.1} satisfies
  \begin{align*}
 \textbf{u}\in C([0,T],\cV)\cap H^2(0,T,L^2(\Omega)^3)\cap H^1(0,T;H^{2}(\Omega)^3),\\
  \textbf{B}\in C([0,T],\cW)\cap H^2(0,T,L^2(\Omega)^3)\cap H^1(0,T;H^{2}(\Omega)^3),\\
 \textbf{w}\in C([0,T],\cV)\cap H^2(0,T,L^2(\Omega)^3)\cap H^1(0,T;H^{2}(\Omega)^3),\\
   p \in C([0,T],\cM\cap H^{1}(\Omega)).
 \end{align*}

We need to introduce Stokes, Maxwell and generalized Ritz projections.
Assume that $(\textbf{u},p)\in \cV\times \cM$, $\textbf{B} \in \cW$, $\textbf{w} \in \cV$,
  we introduce Stokes projection/stabilization Stokes projections $(P_h\textbf{u},J_hp)$ and $(P^{'}_h\textbf{u},J^{'}_hp)$:
\begin{align}\label{eq17}
a_f(\textbf{u}-P_h\textbf{u},\textbf{v})+d(\textbf{v},p-J_hp)&=0, \ \ \forall\ \textbf{v}\in  \cV^1_h,\\\label{eq18}
d(\textbf{u}-P_h\textbf{u},q)&=0, \ \ \forall\ q\in  \cM_h,\\\label{eq19}
\overline{\mathbb{B}}(\textbf{u}-P^{'}_h\textbf{u},p-J^{'}_hp;\textbf{v},q)&=0, \ \ \forall\ (\textbf{v},q)\in  (\cV^2_h,\cM_h).
 \end{align}
 Applying duality argument and Agmon's inequality, we have the approximation properties as follows:
\begin{align}\label{eq20}
\|\nabla(\textbf{u}-Z_h\textbf{u})\|+\|p-Y_hp\|&\leq Ch\bigl(\|\textbf{u}\|_{2}+\|p\|_1\bigr),\\\label{eq21}
\|Z_h\textbf{u}\|_{L^{\infty}}+\|\nabla Z_h\textbf{u}\|_{L^3}&\leq C\bigl(\|\textbf{u}\|_{2}+\|p\|_1\bigr),
 \end{align}
 where $(Z_h\textbf{u},Y_hp)=(P_h\textbf{u},J_hp)$ or $(P^{'}_h\textbf{u},J^{'}_hp)$.
Similarly, we introduce the Maxwell projection $P_{mh}\textbf{B}$:
\begin{align}\label{eq22}
a_B(\textbf{B}-P_{mh}\textbf{B},\textbf{H})&=0, \ &&\forall \ \textbf{H}\in  \cW^1_h\,\mbox{or}\,\cW^2_h,\\\label{eq23}
\|\nabla(\textbf{B}-P_{mh}\textbf{B})\|&\leq Ch\|\textbf{B}\|_{2},\\\label{eq24}
\|P_{mh}\textbf{B}\|_{L^{\infty}}+\|\nabla P_{mh}\textbf{B}_h\|_{L^3}&\leq C\|\textbf{B}\|_{2},
 \end{align}
and the generalized Ritz projection $r_{h}\textbf{w}$:
\begin{align}\label{eq25}
a_w(\textbf{w}-r_h\textbf{w},\psi)&=0, \ &&\forall \ {\psi}\in \cQ^1_h\,\mbox{or}\,\cQ^2_h,\\\label{eq26}
\|\nabla(\textbf{w}-r_h\textbf{w})\|&\leq Ch\|\textbf{w}\|_{2},\\\label{eq27}
\|r_h\textbf{w}\|_{L^{\infty}}+\|\nabla r_h\textbf{w}\|_{L^3}&\leq C\|\textbf{w}\|_{2}.
 \end{align}

 In order to carry out the error analysis, we denote the following discrete errors
\begin{align*}
e^n_\textbf{u} :&= \textbf{u}^n - \textbf{u}^n_h,\ &&e^n_p := p^n -p^n_h,\\
e^n_\textbf{B} :&= \textbf{B}^n- \textbf{B}^n_h, \ &&e^n_\textbf{w} := \textbf{w}^n -\textbf{w}^n_h.
 \end{align*}
The above error estimates can be defined as
\begin{align*}
 e^n_\textbf{u}&=\zeta^n_\textbf{u}-\epsilon^n_\textbf{u},  &&e^n_p=\zeta^n_p-\xi^n_p,\\
  e^n_\textbf{B}&=\zeta^n_\textbf{B}-\xi^n_\textbf{B},  &&e^n_\textbf{w}=\zeta^n_\textbf{w}-\xi^n_\textbf{w},
 \end{align*}
 with
 \begin{align*}
 &\zeta^n_\textbf{u}:=\textbf{u}^n-P_h\textbf{u}^n \in \cV,\ \ \    && \xi_\textbf{u}^n:=\textbf{u}_h^n-P_h\textbf{u}^n \in \cV^1_h\ \mbox{or}\,\ \cV^2_h\\
&\zeta^n_p:=p^n-J_hp^n \in \cM,                     &&\xi_p^n:=p_h^n-J_hp^n \in \cM_h,\\
 &\zeta^n_\textbf{B}:=\textbf{B}^n-P_{mh}\textbf{B}^n \in \cW,           &&\xi_\textbf{B}^n:= \textbf{B}^n_h-P_{mh}\textbf{B}^n \in \cW^1_h\,\ \mbox{or}\,\ \cW^2_h,\\
&\zeta^n_\textbf{w}:=\textbf{w}^n-r_h\textbf{w}^n \in \cV, &&\xi_\textbf{w}^n:=\textbf{w}^n_h-r_h\textbf{w}^n\in \cQ^1_h\,\ \mbox{or}\,\ \cQ^2_h.
 \end{align*}

\begin{theorem}   Suppose that \textbf{Assumpation A2} holds and the initial conditions satisfy
 \begin{align}\label{eq28}
\|\textbf{u}^0-\textbf{u}^0_h\|+ \|\textbf{B}^0-\textbf{B}^0_h\|+\|\textbf{w}^0-\textbf{w}^0_h\|\leq Ch.
 \end{align}
Then the error estimates hold
\begin{align*}
\max_{1\leq n\leq N}\|\textbf{u}(t_n)-\textbf{u}^n_h\|^2
+\Delta t\sum_{n=1}^{N}\|\nabla (\textbf{u}(t_n)-\textbf{u}^n_h)\|^2 &\leq C\bigl(\Delta t^2+h^{2}\bigr),\\
\max_{1\leq n\leq N}\|\textbf{B}(t_n)-\textbf{B}^n_h\|^2
+\Delta t\sum_{n=1}^{N}\|\nabla (\textbf{B}(t_n)-\textbf{B}^n_h)\|^2 &\leq C\bigl(\Delta t^2+h^{2}\bigr),\\
\max_{1\leq n\leq N}\|\textbf{w}(t_n)-\textbf{w}^n_h\|^2
+\Delta t\sum_{n=1}^{N}\|\nabla (\textbf{w}(t_n)-\textbf{w}^n_h)\|^2 &\leq C\bigl(\Delta t^2+h^{2}\bigr).
 \end{align*}
\end{theorem}
\noindent \textit{Proof:}\quad Utilizing \eqref{eq9.1}, \eqref{eq11.1},  \eqref{eq17}-\eqref{eq18}, \eqref{eq22} and \eqref{eq25},  we get the following error equations:
 \begin{subequations}\label{eq29.1}
\begin{alignat}{2} \label{eq29.1a} (\mathcal{D}\xi_\textbf{u}^n,\textbf{v}_h)+a_f(\xi_\textbf{u}^n,\textbf{v}_h)
-d(\textbf{v}_h,\xi_p^n)&=(\Theta_h^n,\textbf{v}_h),\\
d(\xi_\textbf{u}^n,q_h)&=0,\\
S(\mathcal{D}\xi_\textbf{B}^n,\textbf{H}_h)+Sa_B(\xi_\textbf{B}^n,\textbf{H}_h)&=(\widehat{\Theta}_h^n,\textbf{H}_h), \\
(\mathcal{D}\xi_\textbf{w}^n,\psi_h)+a_w(\xi_\textbf{w}^n,\psi_h)+2\nu_r(\xi_\textbf{w}^n,\psi_h)&=(\widetilde{\Theta}_h^n,\psi_h),
\end{alignat}
\end{subequations}
where $(\Theta_h^n,\textbf{v}_h)$, $(\widehat{\Theta}_h^n,\textbf{H}_h)$, $(\widetilde{\Theta}_h^n,\psi_h)$ are denoted by
\begin{align}\label{eq30}
(\Theta_h^n,\textbf{v}_h): &=\bigl(\partial_t\textbf{u}^n-\mathcal{D}P_h{\textbf{u}}^n,\textbf{v}_h\bigr)\\&\nonumber\quad+
e(\textbf{w}^{n}_h,\textbf{v}_h)-e(\textbf{w}^n,\textbf{v}_h)\\&\nonumber\quad
+b(\textbf{u}^n,\textbf{u}^n,\textbf{v}_h)-b(\textbf{u}^{n-1}_h,\textbf{u}^n_h,\textbf{v}_h)\\&\nonumber
  \quad+Sc_{\widehat{B}}(\textbf{B}^n, \textbf{B}^n,\textbf{v}_h)-Sc_{\widehat{B}}(\textbf{B}^{n-1}_h, \textbf{B}^n_h,\textbf{v}_h),
 \end{align}
\begin{align}\label{eq31}
(\widehat{\Theta}_h^n,\textbf{H}_h): &=S\bigl(\partial_t\textbf{B}^n-\mathcal{D}P_{mh}{\textbf{B}}^n,\textbf{H}_h\bigr)\\&\nonumber\quad
  \quad+Sc_{\widehat{B}}(\textbf{u}^{n}_h, \textbf{B}^{n-1}_h,\textbf{H}_h)-Sc_{\widehat{B}}(\textbf{u}^n, \textbf{B}^n,\textbf{v}_h),
 \end{align}
 and
 \begin{align}\label{eq32}
(\widetilde{\Theta}_h^n,\psi_h): &=\bigl(\partial_t\textbf{w}^n-\mathcal{D}r_h\textbf{w}^n,\psi_h\bigr)\\&\nonumber\quad
+e(\textbf{u}_h^n,\psi_h)-e(\textbf{u}^n, \psi_h)+2\nu_r(\zeta_\textbf{w}^n,\psi_h)\\&\nonumber\quad
+c_w(\textbf{u}^n, \textbf{w}^n,\psi_h)-c_w(\textbf{u}^{n-1}_h, \textbf{w}^{n}_h,\psi_h).
 \end{align}
Setting $(\textbf{v}_h, q_h, \textbf{H}_h,\psi_h)=(\xi_\textbf{u}^n, \xi_p^n, \xi_\textbf{B}^n,\xi_\textbf{w}^n)$ in \eqref{eq29.1}, it follows that
 \begin{subequations}\label{eq33.1}
\begin{alignat}{2} \label{eq33.1a}
 (\mathcal{D}\xi_\textbf{u}^n,\xi_\textbf{u}^n)+(\nu+\nu_r)\|\nabla \xi_\textbf{u}^n\|^2&\leq(\Theta_h^n,\xi_\textbf{u}^n),\\
S(\mathcal{D}\xi_\textbf{B}^n,\xi_\textbf{B}^n)+c_{\star}\mu S\|\nabla \xi_\textbf{B}^n\|^2&\leq(\widehat{\Theta}_h^n,\xi_\textbf{B}^n), \\
(\mathcal{D}\xi_\textbf{w}^n,\xi_\textbf{w}^n)+(c_a+c_d)\|\nabla \xi_\textbf{w}^n\|^2&+(c_0+c_d-c_a)\|\nabla \cdot\xi_\textbf{w}^n\|^2\\&\nonumber+2\nu_r\|\xi_\textbf{w}^n\|^2\leq(\widetilde{\Theta}_h^n,\xi_\textbf{w}^n).
\end{alignat}
\end{subequations}
   For the term $(\Theta_h^n,\xi_\textbf{u}^n)$, we rewrite as follows:
 \begin{align}\label{eq34}
 (\Theta_h^n,\xi_\textbf{u}^n)&=\bigl(\partial_t\textbf{u}^n-\mathcal{D}P_h{\textbf{u}}^n,\xi_\textbf{u}^n\bigr)
+e(\xi_\textbf{w}^n,\xi_\textbf{u}^n)+e(\zeta_\textbf{w}^n,\xi_\textbf{u}^n)\\&\nonumber\quad
+b(\textbf{u}^n,\textbf{u}^n-P_h\textbf{u}^n,\xi_\textbf{u}^n)+b(\textbf{u}^n-\textbf{u}^{n-1},P_h\textbf{u}^n,\xi_\textbf{u}^n)
\\&\nonumber
 \quad+b(\textbf{u}^{n-1}-P_h\textbf{u}^{n-1},P_h\textbf{u}^n,\xi_\textbf{u}^n)-b(\xi_\textbf{u}^{n-1},P_h\textbf{u}^n,\xi_\textbf{u}^n)\\&\nonumber
 \quad+Sc_{\widehat{B}}(\textbf{B}^n,\textbf{B}^n-P_{mh}\textbf{B}^n,\xi_\textbf{u}^n)
 +Sc_{\widehat{B}}(\textbf{B}^n-\textbf{B}^{n-1},P_{mh}\textbf{B}^n,\xi_\textbf{u}^n)\\&\nonumber
  \quad+Sc_{\widehat{B}}(\textbf{B}^{n-1}-P_{mh}\textbf{B}^n,P_{mh}\textbf{B}^n,\xi_\textbf{u}^n)-Sc_{\widehat{B}}(\xi_\textbf{B}^{n-1},P_{mh}\textbf{B}^n,\xi_\textbf{u}^n)
   \\&\nonumber
   \quad- b(P_{h}\textbf{u}^{n-1}, \xi_\textbf{u}^n,\xi_\textbf{u}^n)-b(\xi_\textbf{u}^{n-1}, \xi_\textbf{u}^n,\xi_\textbf{u}^n)\\&\nonumber\quad
  - Sc_{\widehat{B}}(P_{mh}\textbf{B}^{n-1}, \xi_\textbf{B}^n,\xi_\textbf{u}^n)
  -Sc_{\widehat{B}}(\xi_\textbf{B}^{n-1}, \xi_\textbf{B}^n,\xi_\textbf{u}^n)\\&\nonumber
  =\sum^{11}_{i=1}\Theta_i -Sc_{\widehat{B}}(P_{mh}\textbf{B}^{n-1}, \xi_\textbf{B}^n,\xi_\textbf{u}^n)-Sc_{\widehat{B}}(\xi_\textbf{B}^{n-1}, \xi_\textbf{B}^n,\xi_\textbf{u}^n).
 \end{align}
 We now bound each term on the RHS of \eqref{eq34}. Using  Young's inequality, \textbf{Assumpation A2} and Taylor's
 theorem, we arrive at
\begin{align*}
\Theta_{1} &\leq C\|\partial_t\textbf{u}^n-\mathcal{D}P_h{\textbf{u}}^n\|\|\nabla \xi_\textbf{u}^n\|\\
&\leq C(h+\Delta t^{\frac{1}{2}}+\frac{h}{\sqrt{\Delta t}})\|{\textbf{u}}\|_{H^1((t_{n-1},t_n],H^1(\Omega)^3)}\|\nabla \xi_\textbf{u}^n\|.
 \end{align*}
Due to \eqref{eq4.1}, Young's inequality, H\"{o}lder inequality and  \eqref{eq21}, we bound the terms $\Theta_2-\Theta_7$ as follows
\begin{align*}
\Bigl|\sum^{7}_{i=2}\Theta_{i}\Bigr| &\leq C\|\zeta_\textbf{w}^n\|\|\nabla\xi_\textbf{u}^n\|+ C\|\xi_\textbf{w}^n\|\|\nabla\xi_\textbf{u}^n\|\\&\nonumber\quad+
C\|\nabla \textbf{u}^n\|\|\nabla\bigl(\textbf{u}^n-P_h\textbf{u}^n\bigr)\|\|\nabla\xi_\textbf{u}^n\|\\&\nonumber\quad
 + C\Delta t^{\frac{1}{2}}\|{\textbf{u}}\|_{H^1((t_{n-1},t_n],H^1(\Omega)^3)}\bigl(\|\nabla P_h\textbf{u}^n\|_{L^3}+\|P_h\textbf{u}^n\|_{L^\infty}\bigr)\|\nabla\xi_\textbf{u}^n\|\\&\nonumber\quad
+C\|\textbf{u}^{n-1}-P_h\textbf{u}^{n-1}\|\bigl(\|\nabla P_h\textbf{u}^n\|_{L^3}+\|P_h\textbf{u}^n\|_{L^\infty}\bigr)\|\nabla \xi_\textbf{u}^n\|\\&\nonumber\quad
+ C\|\xi_\textbf{u}^{n-1}\|\bigl(\|\nabla P_h\textbf{u}^n\|_{L^3}+\|P_h\textbf{u}^n\|_{L^\infty}\bigr)\|\nabla \xi_\textbf{u}^n\|.
\end{align*}
Using \eqref{eq4.1},  H\"{o}lder inequality, Young's inequality and  \eqref{eq24},  the terms $\Theta_8-\Theta_{11}$ can be estimated by
  \begin{align*}
\Bigl|\sum^{11}_{i=8}\Theta_{i}\Bigr|&\leq C\|\textbf{B}^n\|_{L^\infty}\|\nabla\bigl(\textbf{B}^n-P_{mh}\textbf{B}^n\bigr)\|\|\nabla \xi_\textbf{u}^n\|\\&\nonumber\quad
%+ C\|\textbf{B}^n-\textbf{B}^{n-1}\|
+C\Delta t^{\frac{1}{2}}\|{\textbf{B}}\|_{H^1((t_{n-1},t_n],H^1(\Omega)^3)}\|\curl P_{mh}\textbf{B}^n\|_{L^\infty}\|\nabla \xi_\textbf{u}^n\|\\&\nonumber\quad
+C\|\textbf{B}^{n-1}-P_{mh}\textbf{B}^{n-1}\|\|\curl P_{mh}\textbf{B}^n\|_{L^\infty}\|\nabla \xi_\textbf{u}^n\|\\&\nonumber\quad
+ C\|\xi_\textbf{B}^{n-1}\|\|\nabla P_{mh}\textbf{B}^n\|_{L^3}\|\nabla \xi_\textbf{u}^n\|.
\end{align*}
Then inserting \eqref{eq34} with above estimates, \eqref{eq20}, \eqref{eq23} and \eqref{eq26}, we arrive at
  \begin{align}\label{eq35}
 (\Theta_h^n,\xi_\textbf{u}^n)&\leq  C\bigl((h+\Delta t^{\frac{1}{2}}+\frac{h}{\sqrt{\Delta t}})(\|{\textbf{u}}\|_{H^1((t_{n-1},t_n],H^1(\Omega)^3)}+\|{\textbf{B}}\|_{H^1((t_{n-1},t_n],H^1(\Omega)^3)})\\\nonumber&\quad
 +\|\xi_{\textbf{u}}^{n-1}\|+\|\xi_{\textbf{B}}^{n-1}\|+\|\xi_{\textbf{w}}^{n}\|\bigr)\|\nabla \xi_\textbf{u}^n\|\\\nonumber&
  \quad- Sc_{\widehat{B}}(P_{mh}\textbf{B}^{n-1}, \xi_\textbf{B}^n,\xi_\textbf{u}^n)-Sc_{\widehat{B}}(\xi_\textbf{B}^{n-1}, \xi_\textbf{B}^n,\xi_\textbf{u}^n)\\\nonumber&
  \leq C\bigl((h^2+\Delta t+\frac{h^2}{\Delta t})(\|{\textbf{u}}\|^2_{H^1((t_{n-1},t_n],H^1(\Omega)^3)}+\|{\textbf{B}}\|^2_{H^1((t_{n-1},t_n],H^1(\Omega)^3)})\\\nonumber&\quad
    +\|\xi_{\textbf{u}}^{n-1}\|^2+\|\xi_{\textbf{B}}^{n-1}\|^2+\|\xi_{\textbf{w}}^{n}\|^2\bigr)+\frac{(\nu+\nu_r)}{2}\|\nabla \xi_\textbf{u}^n\|^2\\\nonumber&
  \quad- Sc_{\widehat{B}}(P_{mh}\textbf{B}^{n-1}, \xi_\textbf{B}^n,\xi_\textbf{u}^n)-Sc_{\widehat{B}}(\xi_\textbf{B}^{n-1}, \xi_\textbf{B}^n,\xi_\textbf{u}^n).
\end{align}
For the term $(\widehat{\Theta}_h^n,\xi_\textbf{B}^n)$, we rewrite as follows:
 \begin{align}\label{eq36}
 (\widehat{\Theta}_h^n,\xi_\textbf{B}^n)&=S\bigl(\partial_t\textbf{B}^n-\mathcal{D}P_{mh}{\textbf{B}}^n,\xi_\textbf{B}^n\bigr)
 -Sc_{\widetilde{B}}(\textbf{u}^{n}-P_h\textbf{u}^{n},\textbf{B}^{n},\xi_\textbf{B}^n)\\&\nonumber
 \quad-Sc_{\widetilde{B}}(P_h\textbf{u}^{n},\textbf{B}^{n}-\textbf{B}^{n-1},\xi_\textbf{B}^n)
 -Sc_{\widetilde{B}}(P_h\textbf{u}^{n},\textbf{B}^{n-1}-P_{mh}\textbf{B}^{n-1},\xi_\textbf{B}^n)\\&\nonumber
 \quad+Sc_{\widetilde{B}}(P_h\textbf{u}^{n},\xi_\textbf{B}^{n-1},\xi_\textbf{B}^n)
+Sc_{\widetilde{B}}(\xi_\textbf{u}^{n},\xi_\textbf{B}^{n-1},\xi_\textbf{B}^n)
 +Sc_{\widetilde{B}}(\xi_\textbf{u}^{n},P_{mh}\textbf{B}^{n-1},\xi_\textbf{B}^n)\\&\nonumber
 =\sum^{5}_{i=1}\widehat{\Theta}_i^n+Sc_{\widetilde{B}}(\xi_\textbf{u}^{n},\xi_\textbf{B}^{n-1},\xi_\textbf{B}^n)
 +Sc_{\widetilde{B}}(\xi_\textbf{u}^{n},P_{mh}\textbf{B}^{n-1},\xi_\textbf{B}^n).
 \end{align}
 Thanks to \eqref{eq4.1}, \textbf{Assumpation A2}, Taylor's
 theorem and \eqref{eq24},  we bound the terms $\widehat{\Theta}_2^n-\widehat{\Theta}_5^n$ as follows:
  \begin{align*}
  \Bigl|\sum^{5}_{i=1}\widehat{\Theta}_{i}\Bigr|&\leq C(h+\Delta t^{\frac{1}{2}}+\frac{h}{\sqrt{\Delta t}})\|{\textbf{B}}\|_{H^1((t_{n-1},t_n],H^1(\Omega)^3)}\|\nabla \xi_\textbf{B}^n\|\\&\nonumber\quad
+ C\|\textbf{u}^n-P_h\textbf{u}^n\|\|\textbf{B}^n\|_{L^{\infty}}\|\nabla \xi_\textbf{B}^n\|\\&\nonumber\quad
+ C\Delta t^{\frac{1}{2}}\|{\textbf{B}}\|_{H^1((t_{n-1},t_n],H^1(\Omega)^3)}\|P_h\textbf{u}^n\|_{L^{\infty}}\|\nabla \xi_\textbf{B}^n\|\\&\nonumber\quad
+ C\|P_h\textbf{u}^n\|_{L^{\infty}}\|\textbf{B}^{n-1}-P_{mh}\textbf{B}^{n-1}\|\|\nabla \xi_\textbf{B}^n\|\\&\nonumber\quad
+ C\|P_h\textbf{u}^n\|_{L^{\infty}}\|\xi_\textbf{B}^{n-1}\|\|\nabla \xi_\textbf{B}^n\|.
\end{align*}
Inserting \eqref{eq36} with above bounds, \eqref{eq20} and \eqref{eq23},  we obtain
   \begin{align}\label{eq37}
 (\widehat{\Theta}_h^n,\xi_\textbf{B}^n)&\leq  C\bigl((h+\Delta t^{\frac{1}{2}}+\frac{h}{\sqrt{\Delta t}})\|{\textbf{B}}\|_{H^1((t_{n-1},t_n],H^1(\Omega)^3)}+\|\xi_{\textbf{B}}^{n-1}\|^2\bigr)\\\nonumber&
  \quad+Sc_{\widetilde{B}}(\xi_\textbf{u}^{n},\xi_\textbf{B}^{n-1},\xi_\textbf{B}^n)
 +Sc_{\widetilde{B}}(\xi_\textbf{u}^{n},P_{mh}\textbf{B}^{n-1},\xi_\textbf{B}^n)\\\nonumber&
 \leq  C\bigl((h^2+\Delta t+\frac{h^2}{\Delta t})\|{\textbf{B}}\|^2_{H^1((t_{n-1},t_n],H^1(\Omega)^3)}+\|\xi_{\textbf{B}}^{n-1}\|^2\bigr)+\frac{c_{\star}\mu S}{2}\|\nabla \xi_\textbf{B}^n\|^2\\\nonumber&
  \quad+Sc_{\widetilde{B}}(\xi_\textbf{u}^{n},\xi_\textbf{B}^{n-1},\xi_\textbf{B}^n)
 +Sc_{\widetilde{B}}(\xi_\textbf{u}^{n},P_{mh}\textbf{B}^{n-1},\xi_\textbf{B}^n).
\end{align}
For the term $(\widetilde{\Theta}_h^n,\xi_\textbf{w}^n)$, we rewrite as follows:
    \begin{align}\label{eq38}
 (\widetilde{\Theta}_h^n,\xi_\textbf{w}^n)&= \bigl(\partial_t\textbf{w}^n-\mathcal{D}r_h\textbf{w}^n,\xi_\textbf{w}^n\bigr)+e(\xi_\textbf{u}^n,\xi_\textbf{w}^n)+e(\zeta_\textbf{u}^n,\xi_\textbf{w}^n)\\&\nonumber\quad
 +c_w(\textbf{u}^{n},\textbf{w}^{n}-P_h\textbf{w}^{n},\xi_\textbf{w}^n)+c_w(\textbf{u}^{n}-\textbf{u}^{n-1},P_h\textbf{w}^{n},\xi_\textbf{w}^n)\\&\nonumber\quad
 +c_w(\textbf{u}^{n-1}-P_h\textbf{u}^{n-1},P_h\textbf{w}^{n},\xi_\textbf{w}^n)-c_w(\xi_\textbf{u}^{n-1},P_h\textbf{w}^{n},\xi_\textbf{w}^n)\\&\nonumber
 =\sum^{7}_{i=1}\widetilde{\Theta}_i.
 \end{align}
   By using \eqref{eq4.1}, \eqref{eq6},  H\"{o}lder inequality,  Young's inequality and  \eqref{eq27}, we have
 \begin{align*}
\Bigl|\sum^{5}_{i=1}\widetilde{\Theta}_{i}\Bigr| &\leq C(h+\Delta t^{\frac{1}{2}}+\frac{h}{\sqrt{\Delta t}})\|{\textbf{w}}\|_{H^1((t_{n-1},t_n],H^1(\Omega)^3)}\|\nabla\xi_\textbf{w}^n\|+2\nu_r\|\zeta_\textbf{w}^n\|\|\nabla\xi_\textbf{w}^n\|\\&\nonumber\quad
+C\|\zeta_\textbf{u}^n\|\|\nabla\xi_\textbf{w}^n\|+ C\|\xi_\textbf{u}^n\|\|\nabla\xi_\textbf{w}^n\|\\&\nonumber\quad+
C\|\nabla \textbf{u}^n\|\|\nabla\bigl(\textbf{w}^n-r_h\textbf{w}^n\bigr)\|\|\nabla\xi_\textbf{w}^n\|\\&\nonumber\quad
 + C\Delta t^{\frac{1}{2}}\|{\textbf{u}}\|_{H^1((t_{n-1},t_n],H^1(\Omega)^3)}\bigl(\|\nabla r_h\textbf{w}^n\|_{L^3}+\|r_h\textbf{w}^n\|_{L^\infty}\bigr)\|\nabla\xi_\textbf{w}^n\|\\&\nonumber\quad
+C\|\textbf{u}^{n-1}-P_h\textbf{u}^{n-1}\|\bigl(\|\nabla r_h\textbf{w}^n\|_{L^3}+\|r_h\textbf{w}^n\|_{L^\infty}\bigr)\|\nabla \xi_\textbf{w}^n\|\\&\nonumber\quad
+ C\|\xi_\textbf{u}^{n-1}\|\bigl(\|\nabla r_h\textbf{w}^n\|_{L^3}+\|r_h\textbf{w}^n\|_{L^\infty}\bigr)\|\nabla \xi_\textbf{w}^n\|.
\end{align*}
Due to \eqref{eq20} and \eqref{eq26}, we derive
  \begin{align}\label{eq39}
 (\widetilde{\Theta}_h^n,\xi_\textbf{w}^n)&\leq  C\bigl((h+\Delta t^{\frac{1}{2}}+\frac{h}{\sqrt{\Delta t}})(\|{\textbf{w}}\|_{H^1((t_{n-1},t_n],H^1(\Omega)^3)}+\|{\textbf{u}}\|_{H^1((t_{n-1},t_n],H^1(\Omega)^3)})\\\nonumber&\quad
 +\|\xi_{\textbf{u}}^{n-1}\|+\|\xi_{\textbf{u}}^{n}\|\bigr)\|\nabla \xi_\textbf{w}^n\|\\\nonumber
 &\leq  C\bigl((h^2+\Delta t+\frac{h^2}{\Delta t})(\|{\textbf{w}}\|^2_{H^1((t_{n-1},t_n],H^1(\Omega)^3)}+\|{\textbf{u}}\|^2_{H^1((t_{n-1},t_n],H^1(\Omega)^3)})\\\nonumber&\quad
 +\|\xi_{\textbf{u}}^{n-1}\|^2+\|\xi_{\textbf{u}}^{n}\|^2\bigr)+\frac{(c_a+c_d)}{2}\|\nabla \xi_\textbf{w}^n\|^2.
\end{align}
  Combining \eqref{eq33.1} with \eqref{eq7}, \eqref{eq35}, \eqref{eq37} and \eqref{eq39}, we arrive at
  \begin{align}\label{eq40}
 (\mathcal{D}\xi_\textbf{u}^n,\xi_\textbf{u}^n)&+S(\mathcal{D}\xi_\textbf{B}^n,\xi_\textbf{B}^n)
 +(\mathcal{D}\xi_\textbf{w}^n,\xi_\textbf{w}^n)\\\nonumber&+\frac{(\nu+\nu_r)}{2}\|\nabla \xi_\textbf{u}^n\|^2+\frac{c_{\star}\mu S}{2}\|\nabla \xi_\textbf{B}^n\|^2+\frac{(c_a+c_d)}{2}\|\nabla \xi_\textbf{w}^n\|^2\\&\nonumber
 +2\nu_r\|\xi_\textbf{w}^n\|^2+(c_0+c_d-c_a)\|\nabla \cdot\xi_\textbf{w}^n\|^2
 \\\nonumber&\leq C(h^{2}+\Delta t+\frac{h^{2}}{\Delta t})(\|{\textbf{u}}\|^2_{H^1((t_{n-1},t_n],H^1(\Omega)^3)}+\|{\textbf{B}}\|^2_{H^1((t_{n-1},t_n],H^1(\Omega)^3)} \\&\nonumber+\|{\textbf{w}}\|^2_{H^1((t_{n-1},t_n],H^1(\Omega)^3)})
+C\bigl(\|\xi_{\textbf{u}}^{n-1}\|^2+
\|\xi_{\textbf{u}}^{n}\|^2 +\|\xi_{\textbf{B}}^{n-1}\|^2+\|\xi_\textbf{w}^{n}\|^2\bigr).
\end{align}
Summing \eqref{eq40} from $n=1$ to $m$, and using the discrete Gr\"{o}nwall Lemma \cite{Heywood1990,He2015}, we get
  \begin{align}\label{eq41}
 \|\xi_\textbf{u}^m\|^2&+ S\|\xi_\textbf{B}^m\|^2+\|\xi_\textbf{w}^m\|^2+\Delta t(\nu+\nu_r)\sum^{m}_{n=1}\|\nabla \xi_\textbf{u}^n\|^2\\\nonumber&
 +\Delta tc_{\star}\mu S\sum^{m}_{n=1}\|\nabla \xi_\textbf{B}^n\|^2
 +(c_a+c_d)\Delta t\sum^{m}_{n=1}\|\nabla \xi_\textbf{w}^n\|^2+\nu_r\Delta t\sum^{m}_{n=1}\|\xi_\textbf{w}^n\|^2\\&\nonumber
 +(c_0+c_d-c_a)\Delta t\sum^{m}_{n=1}\|\nabla \cdot\xi_\textbf{w}^n\|^2
\leq C_{*}\bigl(h^{2}+\Delta t^2\bigr),
\end{align}
where $C_{*}=C(\|{\textbf{u}}\|^2_{H^1((0,T],H^1(\Omega)^3)},\|{\textbf{B}}\|^2_{H^1((0,T],H^1(\Omega)^3)},\|{\textbf{w}}\|^2_{H^1((0,T],H^1(\Omega)^3)})$.
Utilizing \eqref{eq20}, \eqref{eq23}, \eqref{eq26} and the triangle inequality, the desired result holds.
 The proof is completed.
 $$\eqno\Box$$

 \begin{theorem}
 Suppose that \textbf{Assumpation A2} holds.  Suppose that the initial conditions satisfy \eqref{eq34} and
 \begin{subequations}\label{eq42.1}
\begin{alignat}{2} \label{eq42.1a}
\|\nabla(\textbf{u}^0-\textbf{u}^0_h)\|,\ \|\nabla(\textbf{B}^0-\textbf{B}^0_h)\|,\ \|\nabla(\textbf{w}^0-\textbf{w}^0_h)\|\leq Ch,\\\label{eq42.1b}
d(\textbf{u}^0_h,q)=0,\ \forall\, q\,\in \cM_h.
\end{alignat}
\end{subequations}
Then we have the following error estimates
\begin{align*}
\max_{1\leq n\leq N}\|\nabla(\textbf{u}(t_n)-\textbf{u}^n_h)\|^2
+\Delta t\sum_{n=1}^{N}\|\partial_t\textbf{u}-\mathcal{D}\textbf{u}^n_h\|^2 &\leq C\bigl(\Delta t^2+h^{2}\bigr),\\
\max_{1\leq n\leq N}\|\nabla(\textbf{B}(t_n)-\textbf{B}^n_h)\|^2
+\Delta t\sum_{n=1}^{N}\|\partial_t\textbf{B}-\mathcal{D}\textbf{B}^n_h\|^2 &\leq C\bigl(\Delta t^2+h^{2}\bigr),\\
\max_{1\leq n\leq N}\|\nabla(\textbf{w}(t_n)-\textbf{w}^n_h)\|^2
+\Delta t\sum_{n=1}^{N}\|\partial_t\textbf{w}-\mathcal{D}\textbf{w}^n_h\|^2 &\leq C\bigl(\Delta t^2+h^{2}\bigr).
 \end{align*}
\end{theorem}
\noindent \textit{Proof:}\quad Setting $(\textbf{v}_h, q_h, \textbf{H}_h,\psi_h)=(\mathcal{D}\xi_\textbf{u}^n,0,\mathcal{D}\xi_\textbf{B}^n,\mathcal{D}\xi_\textbf{w}^n)$
in \eqref{eq29.1},  we get
 \begin{subequations}\label{eq43.1}
\begin{alignat}{2} \label{eq43.1a}
 (\mathcal{D}\xi_\textbf{u}^n,\mathcal{D}\xi_\textbf{u}^n)+a_f(\xi_\textbf{u}^n,\mathcal{D}\xi_\textbf{u}^n) &=(\Upsilon_h^n,\mathcal{D}\xi_\textbf{u}^n),\\\label{eq43.1b}
S(\mathcal{D}\xi_\textbf{B}^n,\mathcal{D}\xi_\textbf{B}^n)+Sa_B(\xi_\textbf{B}^n,\mathcal{D}\xi_\textbf{B}^n) &=(\widehat{\Upsilon}_h^n,\mathcal{D}\xi_\textbf{B}^n), \\\label{eq43.1c}
(\mathcal{D}\xi_\textbf{w}^n,\mathcal{D}\xi_\textbf{w}^n)+a_w(\xi_\textbf{w}^n,\mathcal{D}\xi_\textbf{w}^n) +2\nu_r(\xi_\textbf{w}^n,\mathcal{D}\xi_\textbf{w}^n)&=(\widetilde{\Upsilon}_h^n,\mathcal{D}\xi_\textbf{w}^n).
\end{alignat}
\end{subequations}
 The second terms/third term of LHS \eqref{eq43.1} has following equalities
 \begin{align}\label{eq44}
 &a_f(\xi_\textbf{u}^n,\mathcal{D}\xi_\textbf{u}^n)
 =\frac{ \Delta t(\nu+\nu_r)}{2} \|\mathcal{D}\nabla\xi_\textbf{u}^n\|^2 %\\\nonumber& \qquad\qquad
 +\frac{\nu+\nu_r}{2\Delta t}\bigl[\|\nabla\xi_\textbf{u}^n\|^2-\|\nabla\xi_\textbf{u}^{n-1}\|^2\bigr],\\\label{eq45}
 &Sa_B(\xi_\textbf{B}^n,\mathcal{D}\xi_\textbf{B}^n)=
 \frac{ S\Delta t\mu}{2} \|\mathcal{D}\curl\xi_\textbf{B}^n\|^2+\frac{ S\Delta t\mu}{2} \|\mathcal{D}\div\xi_\textbf{B}^n\|^2
 \\\nonumber&\qquad\qquad
 +\frac{S\mu}{2\Delta t}\bigl[(\|\curl\xi_\textbf{B}^n\|+\|\div\xi_\textbf{B}^n\|)-(\|\curl\xi_\textbf{B}^{n-1}\|+\|\div\xi_\textbf{B}^{n-1}\|)\bigr],\\\label{eq46}
 &a_w(\xi_\textbf{w}^n,\mathcal{D}\xi_\textbf{w}^n)=\frac{ \Delta t(c_a+c_d)}{2} \|\mathcal{D}\nabla\xi_\textbf{w}^n\|^2+\frac{ \Delta t(c_0+c_d-c_a)}{2} \|\mathcal{D}\nabla\cdot\xi_\textbf{w}^n\|^2
 \\\nonumber&\qquad\qquad
 +\frac{1}{2\Delta t}\bigl[(c_a+c_d)(\|\nabla\xi_\textbf{w}^n\|-\|\nabla\xi_\textbf{w}^{n-1}\|) \\\nonumber& \qquad\qquad\qquad
 +(c_0+c_d-c_a)(\|\nabla\cdot\xi_\textbf{w}^n\|-\|\nabla\cdot\xi_\textbf{w}^{n-1}\|)\bigr],\\\label{eq47}
 &2\nu_r(\xi_\textbf{w}^n,\mathcal{D}\xi_\textbf{w}^n)
 ={ \Delta t\nu_r} \|\mathcal{D}\xi_\textbf{w}^n\|^2 %\\\nonumber& \qquad\qquad
 +\frac{\nu_r}{\Delta t}\bigl[\|\xi_\textbf{w}^n\|^2-\|\xi_\textbf{w}^{n-1}\|^2\bigr].
 \end{align}
 The first term of RHS in \eqref{eq43.1a} can be rewritten
 \begin{align}\label{eq48}
 &(\Upsilon_h^n,\mathcal{D}\xi_\textbf{u}^n)=\bigl(\partial_t\textbf{u}^n-\mathcal{D}P_h{\textbf{u}}^n,\mathcal{D}\xi_\textbf{u}^n\bigr)
+e(\xi_\textbf{w}^n,\mathcal{D}\xi_\textbf{u}^n)+e(\zeta_\textbf{w}^n,\mathcal{D}\xi_\textbf{u}^n)\\&\nonumber
   \quad+b(\textbf{u}^n,\textbf{u}^n-P_h\textbf{u}^n,\mathcal{D}\xi_\textbf{u}^n)
+b(\textbf{u}^n-\textbf{u}^{n-1},P_h\textbf{u}^n,\mathcal{D}\xi_\textbf{u}^n)\\&\nonumber\quad +b(\textbf{u}^{n-1}-P_h\textbf{u}^{n-1},P_h\textbf{u}^n,\mathcal{D}\xi_\textbf{u}^n)
+b(\xi_\textbf{u}^{n-1},P_h\textbf{u}^n,\mathcal{D}\xi_\textbf{u}^n)\\&\nonumber\quad +Sc_{\widehat{B}}(\textbf{B}^n,\textbf{B}^n-P_{mh}\textbf{B}^n,\mathcal{D}\xi_\textbf{u}^n)
+Sc_{\widehat{B}}(\textbf{B}^n-\textbf{B}^{n-1},P_{mh}\textbf{B}^n,\mathcal{D}\xi_\textbf{u}^n)\\&\nonumber
  \quad+Sc_{\widehat{B}}(\textbf{B}^{n-1}-P_{mh}\textbf{B}^n,P_{mh}\textbf{B}^n,\mathcal{D}\xi_\textbf{u}^n)
  -Sc_{\widehat{B}}(\xi_\textbf{B}^{n-1},P_{mh}\textbf{B}^n,\mathcal{D}\xi_\textbf{u}^n)\\&\nonumber
   \quad- b({\color{blue}P_{h}}\textbf{u}^{n-1}, \xi_\textbf{u}^n,\mathcal{D}\xi_\textbf{u}^n)-b(\xi_\textbf{u}^{n-1}, \xi_\textbf{u}^n,\mathcal{D}\xi_\textbf{u}^n)\\&\nonumber\quad- Sc_{\widehat{B}}(P_{mh}\textbf{B}^{n-1}, \xi_\textbf{B}^n,\mathcal{D}\xi_\textbf{u}^n)
  -Sc_{\widehat{B}}(\xi_\textbf{B}^{n-1}, \xi_\textbf{B}^n,\mathcal{D}\xi_\textbf{u}^n)+d(\mathcal{D}\xi_\textbf{u}^n,\xi_p^n)\\&\nonumber
  =\sum^{16}_{i=1}\Upsilon_i.
 \end{align}
 Using \eqref{eq4.1}, \eqref{eq26} and \eqref{eq27}, the terms $\Upsilon_{1}-\Upsilon_{11}$ can be estimated by
  \begin{align*}
  \Bigl|\sum^{11}_{i=1}\Upsilon_{i}\Bigr|&\leq C\|\mathcal{D}\xi_\textbf{u}^n\|\|\partial_t\textbf{u}^n-\mathcal{D}P_h{\textbf{u}}^n\|\\&\quad
  +C\|\nabla\zeta_\textbf{w}^n\|\|\mathcal{D}\xi_\textbf{u}^n\|+ C\|\nabla\xi_\textbf{w}^n\|\|\mathcal{D}\xi_\textbf{u}^n\|\\&\quad
+ C\|\textbf{u}^n\|_{L^\infty}\|\nabla(\textbf{u}^n-P_h\textbf{u}^n) \|\|\mathcal{D} \xi_\textbf{u}^n\|\\&\quad
+C\|\nabla(\textbf{u}^n-\textbf{u}^{n-1})\|\bigl(\|\nabla P_h\textbf{u}^n\|_{L^3}+\|P_h\textbf{u}^n\|_{L^\infty}\bigr)\|\mathcal{D} \xi_\textbf{u}^n\|\\&\quad
+C\|\nabla(\textbf{u}^{n-1}-P_h\textbf{u}^{n-1})\|\bigl(\|\nabla P_h\textbf{u}^n\|_{L^3}+\|P_h\textbf{u}^n\|_{L^\infty}\bigr)\|\mathcal{D}\xi_\textbf{u}^n\|\\&\quad
+ C\|\nabla\xi_\textbf{u}^{n-1}\|\bigl(\|\nabla P_h\textbf{u}^n\|_{L^3}+\|P_h\textbf{u}^n\|_{L^\infty}\bigr)\|\mathcal{D} \xi_\textbf{u}^n\|\\&\quad
+C\|\textbf{B}^n\|_{L^\infty}\|\nabla(\textbf{B}^n-P_{mh}\textbf{B}^n) \|\|\mathcal{D} \xi_\textbf{u}^n\|\\&\quad
+ C\|\textbf{B}^n-\textbf{B}^{n-1}\|_{L^6}\|\curl P_{mh}\textbf{B}^n\|_{L^3}\|\mathcal{D}\xi_\textbf{u}^n\|\\&\quad
+C\|\textbf{B}^{n-1}-P_{mh}\textbf{B}^{n-1}\|_{L^6}\|\curl P_{mh}\textbf{B}^n\|_{L^3}\|\mathcal{D}\xi_\textbf{u}^n\|\\&\quad
+C\|\xi_\textbf{B}^{n-1}\|_{L^6}\|\curl P_{mh}\textbf{B}^n\|_{L^3}\|\mathcal{D}\xi_\textbf{u}^n\|.
\end{align*}
Applying \textbf{Lemma 3.2}, the terms $\Upsilon_{12}-\Upsilon_{15}$ can be bounded by
  \begin{align*}
 \Bigl|\sum^{15}_{i=12}\Upsilon_{i}\Bigr| &\leq C\|P_h\textbf{u}^{n-1}\|_{L^\infty}\|\nabla \xi_\textbf{u}^n \|\|\mathcal{D} \xi_\textbf{u}^n\|\\& \quad
+C\gamma_{1,n}\|\nabla\xi_\textbf{u}^{n-1}\|\bigl(\|\mathcal{D}\xi_\textbf{u}^n\|+\sqrt{\Delta t}\|\nabla\mathcal{D}\xi_\textbf{u}^n\|\bigr)\\&\quad
+ C\|P_{mh}\textbf{B}^{n-1}\|_{L^\infty}\|\nabla \xi_\textbf{B}^n \|\|\mathcal{D} \xi_\textbf{u}^n\| \\&\quad
+ C\gamma_{2,n}\|\nabla\xi_\textbf{B}^{n-1}\|\bigl(\|\mathcal{D}\xi_\textbf{u}^n\|+\sqrt{\Delta t}\|\nabla\mathcal{D}\xi_\textbf{u}^n\|\bigr),
\end{align*}
where
  \begin{align*}
 \gamma_{1,n}&=\min\bigl\{h^{-1/2},\Delta t^{-1/2}\bigr\}\|\nabla\xi_\textbf{u}^{n}\|,\\
 \gamma_{2,n}&=\min\bigl\{h^{-1/2},\Delta t^{-1/2}\bigr\}\|\nabla\xi_\textbf{B}^{n}\|.
\end{align*}
Due to \eqref{eq42.1b}, we have for $n\geq 1$
  \begin{align*}
\Upsilon_{16}^n=d(\mathcal{D}\xi_\textbf{u}^n,\xi_p^n)=0.
\end{align*}
Combining the above estimates into  \eqref{eq48}, and using \eqref{eq20}, \eqref{eq23}, \eqref{eq26} and Young's inequality, we get
  \begin{align}\label{eq49}
 (\Upsilon_h^n,\mathcal{D}\xi_\textbf{u}^n)&\leq  C_1\bigl(h+\Delta t^{1/2}+\frac{h}{\sqrt{\Delta t}}+\|\nabla\xi_{\textbf{u}}^{n-1}\|+\|\nabla\xi_{\textbf{w}}^{n}\|+\|\nabla\xi_{\textbf{B}}^{n-1}\|\bigr)\|\mathcal{D} \xi_\textbf{u}^n\|
\\\nonumber&
  \quad +C\gamma_{1,n}\|\nabla\xi_\textbf{u}^{n-1}\|\bigl(\|\mathcal{D}\xi_\textbf{u}^n\|+\sqrt{\Delta t}\|\nabla\mathcal{D}\xi_\textbf{u}^n\|\bigr) \\\nonumber& \quad
+ C\gamma_{2,n}\|\nabla\xi_\textbf{B}^{n-1}\|\bigl(\|\mathcal{D}\xi_\textbf{u}^n\|+\sqrt{\Delta t}\|\nabla\mathcal{D}\xi_\textbf{u}^n\|\bigr)\\\nonumber&
  \leq C_1\bigl(h^{2}+\Delta t+\frac{h^{2}}{\Delta t}+\|\nabla\xi_{\textbf{u}}^{n-1}\|^2+\|\nabla\xi_{\textbf{w}}^{n}\|^2+\|\nabla\xi_{\textbf{B}}^{n-1}\|^2  \bigr)\\\nonumber&
  \quad +\frac{1}{8}\|\mathcal{D} \xi_\textbf{u}^n\|^2
+C\gamma_{1,n}\|\nabla\xi_\textbf{u}^{n-1}\|\bigl(\|\mathcal{D}\xi_\textbf{u}^n\|+\sqrt{\Delta t}\|\nabla\mathcal{D}\xi_\textbf{u}^n\|\bigr) \\\nonumber& \quad
+ C\gamma_{2,n}\|\nabla\xi_\textbf{B}^{n-1}\|\bigl(\|\mathcal{D}\xi_\textbf{u}^n\|+\sqrt{\Delta t}\|\nabla\mathcal{D}\xi_\textbf{u}^n\|\bigr)
\\\nonumber&
  \leq C_1\bigl(h^{2}+\Delta t+\frac{h^{2}}{\Delta t}+\|\nabla\xi_{\textbf{u}}^{n-1}\|^2+\|\nabla\xi_{\textbf{w}}^{n}\|^2+\|\nabla\xi_{\textbf{B}}^{n-1}\|^2\bigr)\\\nonumber&
  \quad  +\frac{1}{4}\|\mathcal{D} \xi_\textbf{u}^n\|^2
+C\gamma^2_{1,n}\|\nabla\xi_\textbf{u}^{n-1}\|^2 \\\nonumber& \quad
+ C\gamma^2_{2,n}\|\nabla\xi_\textbf{B}^{n-1}\|^2+\frac{\Delta t(\nu+\nu_r)}{4}\|\nabla\mathcal{D}\xi_\textbf{u}^n\|^2,
\end{align}
where $C_1=C(\|{\textbf{u}}\|_{H^1((t_{n-1},t_n],H^1(\Omega)^3)},\|{\textbf{B}}\|_{H^1((t_{n-1},t_n],H^1(\Omega)^3)})$.
The first term of RHS in  \eqref{eq43.1b} can be rewritten
 \begin{align}\label{eq50}
 &(\widehat{\Upsilon}_h^n,\mathcal{D}\xi_\textbf{B}^n) = \bigl(\partial_t\textbf{B}^n-\mathcal{D}P_{mh}{\textbf{B}}^n,\mathcal{D}\xi_\textbf{B}^n\bigr)\\&\nonumber\quad
 - Sc_{\widetilde{B}}(\textbf{u}^{n}-P_h\textbf{u}^{n},\textbf{B}^{n},\mathcal{D}\xi_\textbf{B}^n)
 -Sc_{\widetilde{B}}(P_h\textbf{u}^{n},\textbf{B}^{n}-\textbf{B}^{n-1},\mathcal{D}\xi_\textbf{B}^n)\\&\nonumber
 \quad-Sc_{\widetilde{B}}(P_h\textbf{u}^{n},\textbf{B}^{n-1}-P_{mh}\textbf{B}^{n-1},\xi_\textbf{B}^n)
 +Sc_{\widetilde{B}}(P_h\textbf{u}^{n},\xi_\textbf{B}^{n-1},\mathcal{D}\xi_\textbf{B}^n)\\&\nonumber
 \quad +Sc_{\widetilde{B}}(\xi_\textbf{u}^{n},P_{mh}\textbf{B}^{n-1},\mathcal{D}\xi_\textbf{B}^n)
 +Sc_{\widetilde{B}}(\xi_\textbf{u}^{n},\xi_\textbf{B}^{n-1},\mathcal{D}\xi_\textbf{B}^n)
\\&\nonumber
 =\sum^{7}_{i=1}\widehat{\Upsilon}_i.
 \end{align}
Thanks to \eqref{eq4.1}, Young's inequality and the H\"{o}lder's inequality,  the terms $\widehat{\Upsilon}_{1}$-$\widehat{\Upsilon}_{5}$ can be bounded by
  \begin{align*}
   \Bigl|\sum^{5}_{i=1}\widehat{\Upsilon}_{i}\Bigr| &\leq C\|\mathcal{D}\xi_\textbf{B}^n\|\|\partial_t\textbf{B}^n-\mathcal{D}R_h{\textbf{B}}^n\|+
S\|\nabla(\textbf{u}^n-P_h\textbf{u}^n)\|\|\textbf{B}^n\|_{L^{\infty}}\|\mathcal{D}\xi_\textbf{B}^n\|\\&\quad
+ S\|\textbf{u}^n-P_h\textbf{u}^n\|_{L^6}\|\nabla\textbf{B}^n\|_{L^{3}}\|\mathcal{D}\xi_\textbf{B}^n\|
 +S\|\nabla P_h\textbf{u}^n\|\|\textbf{B}^n-\textbf{B}^{n-1}\|_{L^{\infty}}\|\mathcal{D}\xi_\textbf{B}^n\|\\&\quad+S\|\nabla P_h\textbf{u}^n\|_{L^3}\|\textbf{B}^n-\textbf{B}^{n-1}\|_{L^{6}}\|\mathcal{D}\xi_\textbf{B}^n\|
+S\|\nabla P_h\textbf{u}^n\|_{L^{3}}\|\textbf{B}^{n-1}-P_{mh}\textbf{B}^{n-1}\|_{L^6}\|\mathcal{D}\xi_\textbf{B}^n\|\\&\quad
+S\|P_h\textbf{u}^n\|_{L^{\infty}}\|\nabla(\textbf{B}^{n-1}-P_{mh}\textbf{B}^{n-1})\|\|\mathcal{D}\xi_\textbf{B}^n\|
+ S\|\nabla P_h\textbf{u}^n\|_{L^{3}}\|\xi_\textbf{B}^{n-1}\|_{L^6}\|\mathcal{D}\xi_\textbf{B}^n\|\\&\quad
+S\|P_h\textbf{u}^n\|_{L^{\infty}}\|\nabla\xi_\textbf{B}^{n-1}\|\|\mathcal{D}\xi_\textbf{B}^n\|.
\end{align*}
 Using  \textbf{Lemma 3.2}, the terms $\widehat{\Upsilon}_{6}$ and $\widehat{\Upsilon}_{7}$ can be estimated by
  \begin{align*}
 \Bigl|\sum^{7}_{i=6}\widehat{\Upsilon}_{i}\Bigr| &\leq S\|\nabla\xi_\textbf{u}^{n}\|\|P_{mh}\textbf{B}^{n-1}\|_{L^\infty}\|\mathcal{D}\xi_\textbf{B}^n\|
  +S\|\xi_\textbf{u}^{n}\|_{L^6}\|\nabla P_{mh}\textbf{B}^{n-1}\|_{L^3}\|\mathcal{D}\xi_\textbf{B}^n\|\\&\quad
  + C\gamma_{3,n}\|\nabla\xi_\textbf{B}^{n-1}\|\bigl(\|\mathcal{D}\xi_\textbf{B}^n\|+\sqrt{\Delta t}\|\nabla\mathcal{D}\xi_\textbf{B}^n\|\bigr),
\end{align*}
where $\gamma_{3,n}=\min\bigl\{h^{-1/2},\Delta t^{-1/2}\bigr\}\|\nabla\xi_\textbf{u}^{n}\|$.

Inserting the above estimates into \eqref{eq50}, and utilizing \eqref{eq23} and Young's inequality, we obtain
   \begin{align}\label{eq51}
 (\widehat{\Upsilon}_h^n,\mathcal{D}\xi_\textbf{B}^n)&\leq  C_2\bigl(h+\Delta t^{1/2}+\frac{h}{\sqrt{\Delta t}}+\|\nabla\xi_{\textbf{B}}^{n-1}\|\bigr)\|\mathcal{D}\xi_\textbf{B}^n\|\\\nonumber&
  \quad+ C\gamma_{3,n}\|\nabla\xi_\textbf{B}^{n-1}\|\bigl(\|\mathcal{D}\xi_\textbf{B}^n\|+\sqrt{\Delta t}\|\nabla\mathcal{D}\xi_\textbf{B}^n\|\bigr)\\\nonumber&
 \leq  C_2\bigl(h^{2}+\Delta t+\frac{h^{2}}{\Delta t}+\|\xi_\textbf{w}^{n-1}\|^2+\|\xi_{\textbf{B}}^{n-1}\|^2+\|\nabla\xi_\textbf{u}^{n}\|^2\bigr)+\frac{1}{8}\|\mathcal{D}\xi_\textbf{B}^n\|^2\\\nonumber&
  \quad+ C\gamma_{3,n}\|\nabla\xi_\textbf{B}^{n-1}\|\bigl(\|\mathcal{D}\xi_\textbf{B}^n\|+\sqrt{\Delta t}\|\nabla\mathcal{D}\xi_\textbf{B}^n\|\bigr) \\\nonumber&
 \leq  C_2\bigl(h^{2}+\Delta t+\frac{h^{2}}{\Delta t}+\|\xi_\textbf{w}^{n-1}\|^2+\|\xi_{\textbf{B}}^{n-1}\|^2+\|\nabla\xi_\textbf{u}^{n}\|^2\bigr)+\frac{1}{4}\|\mathcal{D}\xi_\textbf{B}^n\|^2\\\nonumber&
  \quad+ C\gamma^2_{3,n}\|\nabla\xi_\textbf{B}^{n-1}\|^2+\frac{\Delta t\mu c_{\star}S}{4}\|\nabla\mathcal{D}\xi_\textbf{B}^n\|^2.
\end{align}
where $C_2=C(\|{\textbf{B}}\|_{H^1((t_{n-1},t_n],H^1(\Omega)^3)})$.
The first term of RHS in \eqref{eq43.1c} can be rewritten
    \begin{align}\label{eq52}
 (\widetilde{\Upsilon}_h^n,\mathcal{D}\xi_\textbf{w}^n)&=\bigl(\partial_t\textbf{w}^n-\mathcal{D}r_{h}{\textbf{w}}^n,\mathcal{D}\xi_\textbf{w}^n\bigr)
 +2\nu_r(\zeta_\textbf{w}^n,\mathcal{D}\xi_\textbf{w}^n)\\&\nonumber\quad
 +e(\xi_\textbf{u}^n,\mathcal{D}\xi_\textbf{w}^n)+e(\zeta_\textbf{u}^n,\mathcal{D}\xi_\textbf{w}^n)\\&\nonumber\quad
+ c_w(\textbf{u}^{n},\textbf{w}^{n},\mathcal{D}\xi_\textbf{w}^n)
+c_w(\textbf{u}^{n}-\textbf{u}^{n-1},r_h\textbf{w}^{n},\mathcal{D}\xi_\textbf{w}^n)\\&\nonumber \quad
+c_w(\textbf{u}^{n-1}-P_h\textbf{u}^{n-1},r_h\textbf{w}^{n},\mathcal{D}\xi_\textbf{w}^n)
-c_w(\xi_\textbf{u}^{n-1},P_h\textbf{w}^{n},\mathcal{D}\xi_\textbf{w}^n)\\&\nonumber
 \quad-c_w(P_h\textbf{u}^{n-1},\xi_\textbf{w}^{n},\mathcal{D}\xi_\textbf{w}^n)
 -c_w(\xi_\textbf{u}^{n-1},\xi_\textbf{w}^{n},\mathcal{D}\xi_\textbf{w}^n)\\&\nonumber
 =\sum^{10}_{i=1}\widetilde{\Upsilon}_i.
 \end{align}
Using  \eqref{eq4.1}, Young's inequality and the H\"{o}lder's inequality,  the terms $\widetilde{\Upsilon}_{1}$--$\widetilde{\Upsilon}_{8}$ can be bounded by
   \begin{align*}
  \Bigl|\sum^{8}_{i=1}\widetilde{\Upsilon}_{i}\Bigr|&\leq C\|\mathcal{D}\xi_\textbf{w}^n\|\|\partial_t\textbf{w}^n-\mathcal{D}r_h{\textbf{w}}^n\|\\&\quad
  +C\|\nabla\zeta_\textbf{u}^n\|\|\mathcal{D}\xi_\textbf{w}^n\|+ C\|\nabla\xi_\textbf{u}^n\|\|\mathcal{D}\xi_\textbf{w}^n\|\\&\quad
+ C\|\textbf{u}^n\|_{L^\infty}\|\nabla(\textbf{w}^n-r_h\textbf{w}^n) \|\|\mathcal{D} \xi_\textbf{w}^n\|\\&\quad
+C\|\nabla(\textbf{u}^n-\textbf{u}^{n-1})\|\bigl(\|\nabla r_h\textbf{w}^n\|_{L^3}+\|r_h\textbf{w}^n\|_{L^\infty}\bigr)\|\mathcal{D} \xi_\textbf{w}^n\|\\&\quad
+C\|\nabla(\textbf{u}^{n-1}-P_h\textbf{u}^{n-1})\|\bigl(\|\nabla r_h\textbf{w}^n\|_{L^3}+\|r_h\textbf{w}^n\|_{L^\infty}\bigr)\|\mathcal{D}\xi_\textbf{w}^n\|\\&\quad
+ C\|\nabla\xi_\textbf{u}^{n-1}\|\bigl(\|\nabla r_h\textbf{w}^n\|_{L^3}+\|r_h\textbf{w}^n\|_{L^\infty}\bigr)\|\mathcal{D} \xi_\textbf{w}^n\|.
\end{align*}
From \textbf{Lemma 3.2}, the terms $\widetilde{\Upsilon}_{9}-\widetilde{\Upsilon}_{10}$ can be estimated by
  \begin{align*}
 \Bigl|\sum^{10}_{i=9}\widetilde{\Upsilon}_{i}\Bigr| &\leq C\|P_h\textbf{u}^{n-1}\|_{L^\infty}\|\nabla \xi_\textbf{w}^n \|\|\mathcal{D} \xi_\textbf{w}^n\|\\& \quad
+C\gamma_{4,n}\|\nabla\xi_\textbf{u}^{n-1}\|\bigl(\|\mathcal{D}\xi_\textbf{w}^n\|+\sqrt{\Delta t}\|\nabla\mathcal{D}\xi_\textbf{w}^n\|\bigr),
\end{align*}
where $\gamma_{4,n}=\min\bigl\{h^{-1/2},\Delta t^{-1/2}\bigr\}\|\nabla\xi_\textbf{w}^{n}\|$.

Inserting the above estimates into \eqref{eq52}, and utilizing \eqref{eq20}, \eqref{eq26} and  Young's inequality, we have
  \begin{align}\label{eq53}
 (\widetilde{\Upsilon}_h^n,\mathcal{D}\xi_\textbf{w}^n)&\leq  C_3\bigl(h+\Delta t^{1/2}+\frac{h}{\sqrt{\Delta t}}+\|\nabla\xi_{\textbf{u}}^{n-1}\|+\|\nabla\xi_{\textbf{u}}^{n}\|\bigr)\|\mathcal{D}\xi_\textbf{w}^n\|\\\nonumber&\quad
 + C\gamma_{4,n}\|\nabla\xi_\textbf{u}^{n-1}\|\bigl(\|\mathcal{D}\xi_\textbf{w}^n\|+\sqrt{\Delta t}\|\nabla\mathcal{D}\xi_\textbf{w}^n\|\bigr)\\\nonumber
 &\leq  C_3\bigl(h^{2}+\Delta t+\frac{h^{2}}{\Delta t}+\|\nabla\xi_{\textbf{u}}^{n-1}\|^2+\|\nabla\xi_{\textbf{u}}^{n}\|^2\bigr)+\frac{1}{8}\|\mathcal{D}\xi_\textbf{w}^n\|^2\\\nonumber&\quad
  + C\gamma_{4,n}\|\nabla\xi_\textbf{u}^{n-1}\|\bigl(\|\mathcal{D}\xi_\textbf{w}^n\|+\sqrt{\Delta t}\|\nabla\mathcal{D}\xi_\textbf{w}^n\|\bigr)\\\nonumber
 &\leq  C_3\bigl(h^{2}+\Delta t+\frac{h^{2}}{\Delta t}+\|\nabla\xi_{\textbf{u}}^{n-1}\|^2+\|\nabla\xi_{\textbf{u}}^{n}\|^2\bigr)+\frac{1}{4}\|\mathcal{D}\xi_\textbf{w}^n\|^2\\\nonumber&\quad
  + C\gamma^2_{4,n}\|\nabla\xi_\textbf{u}^{n-1}\|^2+\frac{\Delta t(c_a+c_d)}{4}\|\nabla\mathcal{D}\xi_\textbf{w}^n\|^2,
\end{align}
where $C_3=C(\|{\textbf{u}}\|_{H^1((t_{n-1},t_n],H^1(\Omega)^3)},\|{\textbf{w}}\|_{H^1((t_{n-1},t_n],H^1(\Omega)^3)})$.
 Combining \eqref{eq43.1} with \eqref{eq44}-\eqref{eq47}, \eqref{eq49}, \eqref{eq51} and \eqref{eq53}, it follows that
  \begin{align}\label{eq54}
 &\|\mathcal{D}\xi_\textbf{u}^n\|^2+S\|\mathcal{D}\xi_\textbf{B}^n\|^2+\|\mathcal{D}\xi_\textbf{w}^n\|^2
 +\frac{\nu+\nu_r}{\Delta t}\bigl[\|\nabla\xi_\textbf{u}^n\|^2-\|\nabla\xi_\textbf{u}^{n-1}\|^2\bigr]\\\nonumber&\quad+
 \frac{S\mu}{\Delta t}\bigl[(\|\curl\xi_\textbf{B}^n\|+\|\div\xi_\textbf{B}^n\|)-(\|\curl\xi_\textbf{B}^{n-1}\|+\|\div\xi_\textbf{B}^{n-1}\|)\bigr]\\\nonumber&\quad +
 \frac{1}{2\Delta t}\bigl[(c_a+c_d)(\|\nabla\xi_\textbf{w}^n\|-\|\nabla\xi_\textbf{w}^{n-1}\|) \\\nonumber& \quad
 +(c_0+c_d-c_a)(\|\nabla\cdot\xi_\textbf{w}^n\|-\|\nabla\cdot\xi_\textbf{w}^{n-1}\|)\bigr] +\frac{2\nu_r}{\Delta t}\bigl[\|\xi_\textbf{w}^n\|^2-\|\xi_\textbf{w}^{n-1}\|^2\bigr]
\\\nonumber&\quad+{\Delta t(\nu+\nu_r)}\|\nabla\mathcal{D}\xi_\textbf{u}^n\|^2
+{\Delta t\mu c_{\star}S}\|\nabla\mathcal{D}\xi_\textbf{B}^n\|^2+{\Delta t(c_a+c_d)}\|\nabla\mathcal{D}\xi_\textbf{w}^n\|^2\\\nonumber&\quad
+{\Delta t(c_0+c_d-c_a)}\|\nabla\cdot\mathcal{D}\xi_\textbf{w}^n\|^2 +{ \Delta t\nu_r} \|\mathcal{D}\xi_\textbf{w}^n\|^2%\\\nonumber& \qquad\qquad
\\\nonumber
 &\leq  C\bigl(\gamma^2_{1,n}+\gamma^2_{4,n}\bigr)\|\nabla\xi_{\textbf{u}}^{n-1}\|^2
 +C\bigl(\gamma^2_{2,n}+\gamma^2_{3,n}\bigr)\|\nabla\xi_{\textbf{B}}^{n-1}\|^2+C_0\varsigma_n,
\end{align}
where
  \begin{align}\label{eq55}
\varsigma_n&=h^{2}+\Delta t+\frac{h^{2}}{\Delta t}+\|\nabla\xi_{\textbf{u}}^{n-1}\|^2%+\|\nabla\xi_{\textbf{B}}^{n}\|^2
\\\nonumber&\quad+\|\nabla\xi_{\textbf{u}}^{n}\|^2+\|\nabla\xi_{\textbf{B}}^{n-1}\|^2
 +\|\nabla\xi_{\textbf{w}}^{n}\|^2.
\end{align}
Summing \eqref{eq54} from $n=1$ to $m$,  we obtain
  \begin{align}\label{eq56}
 &\Delta t\sum_{n=1}^m\|\mathcal{D}\xi_\textbf{u}^n\|^2+S\Delta t\sum_{n=1}^m\|\mathcal{D}\xi_\textbf{B}^n\|^2+\Delta t\sum_{n=1}^m\|\mathcal{D}\xi_\textbf{w}^n\|^2 \\\nonumber&\quad
 +(\nu+\nu_r)\|\nabla\xi_\textbf{u}^m\|^2+S\mu c_{\star}\|\nabla\xi_\textbf{B}^m\|^2+(c_a+c_d)\|\nabla\xi_\textbf{w}^m\|^2
 \\\nonumber&\quad+(c_0+c_d-c_a)\|\nabla\cdot\xi_\textbf{w}^m\|^2+2\nu_r\|\xi_\textbf{w}^m\|^2\\\nonumber
 &\leq % C\sum_{n=1}^m\|\nabla\xi^{n-1}_\textbf{u}\|^2
% +C\sum_{n=1}^m\|\nabla\xi^{n-1}_\textbf{B}\|^2+
  C\sum_{n=1}^m\|\nabla\xi^{n-1}_\textbf{w}\|^2\\\nonumber&\quad
 +C\bigl(\|\nabla\xi^{0}_\textbf{u}\|^2+\|\nabla\xi^{0}_\textbf{B}\|^2+\|\nabla\xi^{0}_\textbf{w}\|^2\bigr)
 + C\Delta t\sum_{n=1}^m\bigl(\gamma^2_{1,n}+\gamma^2_{4,n}\bigr)\|\nabla\xi_{\textbf{u}}^{n-1}\|^2\\\nonumber&\quad
 +C\Delta t\sum_{n=1}^m\bigl(\gamma^2_{2,n}+\gamma^2_{3,n}\bigr)\|\nabla\xi_{\textbf{B}}^{n-1}\|^2+\widehat{C}\Delta t\sum_{n=1}^m\varsigma_n,
\end{align}
where $\widehat{C}=C(\|{\textbf{u}}\|_{H^1((0,T],H^1(\Omega)^3)},\|{\textbf{B}}\|_{H^1((0,T],H^1(\Omega)^3)},\|{\textbf{w}}\|_{H^1((0,T],H^1(\Omega)^3)})$.
By using  \textbf{Theorem 4.1}, we have
 \begin{align}\label{eq57}
  &\Delta t\sum_{n=1}^m\|\mathcal{D}\xi_\textbf{u}^n\|^2+S\Delta t\sum_{n=1}^m\|\mathcal{D}\xi_\textbf{B}^n\|^2+\Delta t\sum_{n=1}^m\|\mathcal{D}\xi_\textbf{w}^n\|^2 \\\nonumber&\quad
 +(\nu+\nu_r)\|\nabla\xi_\textbf{u}^m\|^2+S\mu c_{\star}\|\nabla\xi_\textbf{B}^m\|^2+(c_a+c_d)\|\nabla\xi_\textbf{w}^m\|^2
 \\\nonumber&\quad+(c_0+c_d-c_a)\|\nabla\cdot\xi_\textbf{w}^m\|^2+2\nu_r\|\nabla\xi_\textbf{w}^m\|^2\\\nonumber
&\leq  C\bigl(\Delta t^2+h^{2}\bigr)+C\Delta t\sum_{n=1}^m\bigl(\gamma^2_{1,n}+\gamma^2_{4,n}\bigr)\|\nabla\xi_{\textbf{u}}^{n-1}\|^2\\\nonumber&\quad
 +C\Delta t\sum_{n=1}^m\bigl(\gamma^2_{2,n}+\gamma^2_{3,n}\bigr)\|\nabla\xi_{\textbf{B}}^{n-1}\|^2+\widehat{C}\Delta t\sum_{n=1}^m\varsigma_n.
\end{align}
Applying \textbf{Theorem 4.1}, one finds that
  \begin{align}\label{eq58}
 \Delta t\sum_{i=1}^m(\gamma^2_{1,n}+\gamma^2_{4,n})&\leq \min\bigl\{h^{-1},\Delta t^{-1}\bigr\}\Delta t \sum_{i=1}^m\bigl(\|\nabla e_\textbf{u}^n\|^2+\|\nabla e_\textbf{B}^n\|^2\bigr)\\\nonumber&
% \leq  \min\bigl\{h^{-d/3},\Delta t^{-1}\bigr\}\bigl(\Delta t^2+h^{2k+2}\bigr)\\\nonumber&
 \leq C\min\bigl\{h,\Delta t\bigr\}\leq C,\\\label{eq71}
 \Delta t\sum_{i=1}^m\bigl(\gamma^2_{2,n}+\gamma^2_{3,n}\bigr)&\leq C,\\\label{eq72}
 \Delta t\sum_{i=1}^m\varsigma_n&\leq C\bigl(\Delta t^2+h^{2}\bigr).
\end{align}
Using the discrete Gr\"{o}nwall Lemma to \eqref{eq57} with \eqref{eq20}, \eqref{eq23}, \eqref{eq26} and the triangle inequality, the desired result holds.
 The proof is completed.
 $$\eqno\Box$$

\begin{theorem} Suppose that assumptions of Theorems 4.1-4.2 hold.
Then we have following error estimate
\begin{align*}
\Delta t\sum_{n=1}^{N}\|p(t_n)-p^n_h\|^2 &\leq C\bigl(\Delta t^2+h^{2}\bigr).
 \end{align*}
\end{theorem}
\noindent \textit{Proof:}\quad  By using \textbf{Assumption A1} and \eqref{eq29.1a}, we have
\begin{align}\label{eq59}
{\beta}_0\|\xi^n_p\|&\leq \sup_{\textbf{v}_h\in \cV_h, \textbf{v}_h\neq\textbf{0}}\Bigl\{\dfrac{(\mathcal{D}\xi_\textbf{u}^n,\textbf{v}_h)+a_f(\xi_\textbf{u}^n,\textbf{v}_h)}{\|\nabla\textbf{v}_h\|}\Bigr\} \\\nonumber&\quad
+\sup_{\textbf{v}_h\in \cV_h, \textbf{v}_h\neq\textbf{0}}\Bigl\{\dfrac{-(\Theta_h^n,\textbf{v}_h)-(\partial_t\textbf{u}^n-\mathcal{D}P_h{\textbf{u}}^n,\textbf{v}_h)}
{\|\nabla\textbf{v}_h\|}\Bigr\}.
 \end{align}
Similar to estimate \eqref{eq35}, we bound \eqref{eq59} as follows
 \begin{align*}
 \|\xi^n_p\|&\leq C\bigr(\Delta t+h^{2}+\|\mathcal{D}\xi_\textbf{u}^n\|+\|\nabla\xi_\textbf{u}^n\|
+ \|\xi_{\textbf{u}}^{n-1}\|\\\nonumber&\quad+\|\xi_{\textbf{B}}^{n-1}\|+\|\xi_\textbf{w}^{n}\|
+\|\nabla\xi_{\textbf{B}}^{n-1}\|+\|\nabla\xi_{\textbf{B}}^{n}\|\bigr).
 \end{align*}
Thanks to \textbf{Theorems 4.1, 4.2} and the triangle inequality, the desired result holds. The proof is finished.
 $$\eqno\Box$$

  \textbf{Remark 2:}
Using similar techniques as the proof of {Theorem 4.1}-{Theorem 4.3}, {Lemma 3.1} and \eqref{eq19},
 we can obtain a similar error estimate result for Scheme  3.2, here we omit it.

 \section{\label{Sec3} Second-order/decoupled  numerical schemes}

 \subsection{\label{Sec3} Second-order numerical schemes }

In this section, based on Crank-Nicolson discretization
in time and extrapolated treatment of the nonlinear terms, we propose
two second-order numerical schemes of problem \eqref{eq9.1} as follows.

\textbf{Scheme  5.1:}  Given $(\textbf{u}^0_h,\textbf{B}^0_h,\textbf{w}^0_h),(\textbf{u}^1_h,\textbf{B}^1_h,\textbf{w}^1_h) \in \cV^1_h\times \cW^1_h\times \cQ^1_h$,
find $(\textbf{u}^n_h, p^n_h, \textbf{B}^n_h, \textbf{w}^n_h) \in {\cV}^1_h\times \cM_h\times \cW^1_h\times \cQ^1_h$, such that
 \begin{subequations}\label{eq60.1}
\begin{alignat}{2} \label{eq60.1a}
(\mathcal{D}\textbf{u}^{n}_h,\textbf{v}_h)+a_f(\overline{\textbf{u}}^n_h,\textbf{v}_h)
+b(\mathcal{I}(\textbf{u}^{n}_h),\overline{\textbf{u}}^n_h,\textbf{v}_h)\\\nonumber\quad
+Sc_{\widehat{B}}(\mathcal{I}(\textbf{B}^{n}_h),\overline{\textbf{B}}^n_h,\textbf{v}_h)-d(\textbf{v}_h,\overline{p}^n_h)
&=e(\overline{\textbf{w}}_h^n,\textbf{v}_h)+(\textbf{f}^{n-\frac{1}{2}},\textbf{v}_h), \\
d(\overline{\textbf{u}}^n_h,q_h)&=0,\\
(\mathcal{D}\textbf{B}^{n}_h,\textbf{H}_h)+a_B(\overline{\textbf{B}}^{n}_h,\textbf{H}_h)
-c_{\widetilde{B}}(\overline{\textbf{u}}^{n}_h,\mathcal{I}(\textbf{B}^{n}_h),\textbf{H}_h)&=0, \\
(\mathcal{D}\textbf{w}^{n}_h,\psi_h)+a_w(\overline{\textbf{w}}^{n}_h,\psi_h)+c_w(\mathcal{I}(\textbf{u}^{n}_h),\overline{\textbf{w}}^{n}_h,\psi_h)
\\\nonumber\quad+2\nu_r(\overline{\textbf{w}}_h^n,\psi_h)&=e(\overline{\textbf{u}}_h^n,\psi_h)+(\textbf{g}^{n-\frac{1}{2}},\psi_h),
\end{alignat}
\end{subequations}
for all $(\textbf{v}_h,q_h,\textbf{H}_h,\psi_h)\in {\cV}^1_h\times \cM_h\times \cW^1_h\times \cQ^1_h$, where $\overline{\textbf{u}}^n_h=\frac{\textbf{u}^n_h+\textbf{u}^{n-1}_h}{2}$,
 $\overline{\textbf{B}}^n_h=\frac{\textbf{B}^n_h+\textbf{B}^{n-1}_h}{2}$, $\overline{\textbf{w}}^n_h=\frac{\textbf{w}^n_h+\textbf{w}^{n-1}_h}{2}$
 and $\mathcal{I}(\psi^n_h)=\frac{3\psi_h^{n-1}-\psi_h^{n-2}}{2}$.

Using stabilized finite element pair, we propose a new scheme as follows:

\textbf{Scheme 5.2:}  Given $(\textbf{u}^0_h,\textbf{B}^0_h,\textbf{w}^0_h), (\textbf{u}^1_h,\textbf{B}^1_h,\textbf{w}^1_h) \in \cV^2_h\times \cW^2_h\times \cQ^2_h$,
find $(\textbf{u}^n_h, p^n_h, \textbf{B}^n_h, \textbf{w}^n_h) \in {\cV}^2_h\times \cM_h\times \cW^2_h\times \cQ^2_h$, such that
 \begin{subequations}\label{eq61.1}
\begin{alignat}{2} \label{eq61.1a}
(\mathcal{D}\textbf{u}^{n}_h,\textbf{v}_h)+\overline{\mathbb{B}}(\overline{\textbf{u}}^n_h,\overline{p}^n_h;\textbf{v}_h,q_h)
+b(\mathcal{I}(\textbf{u}^{n}_h),\overline{\textbf{u}}^n_h,\textbf{v}_h)\\\nonumber\quad
+Sc_{\widehat{B}}(\mathcal{I}(\textbf{B}^{n}_h),\overline{\textbf{B}}^n_h,\textbf{v}_h)%-d(\textbf{v}_h,\overline{p}^n_h)
&=e(\overline{\textbf{w}}_h^n,\textbf{v}_h)+(\textbf{f}^{n-\frac{1}{2}},\textbf{v}_h), \\
%d(\overline{\textbf{u}}^n_h,q_h)&=0,\\
(\mathcal{D}\textbf{B}^{n}_h,\textbf{H}_h)+a_B(\overline{\textbf{B}}^{n}_h,\textbf{H}_h)
-c_{\widetilde{B}}(\overline{\textbf{u}}^{n}_h,\mathcal{I}(\textbf{B}^{n}_h),\textbf{H}_h)&=0, \\
(\mathcal{D}\textbf{w}^{n}_h,\psi_h)+a_w(\overline{\textbf{w}}^{n}_h,\psi_h)+c_w(\mathcal{I}(\textbf{u}^{n}_h),\overline{\textbf{w}}^{n}_h,\psi_h)
\\\nonumber\quad+2\nu_r(\overline{\textbf{w}}_h^n,\psi_h)&=e(\overline{\textbf{u}}_h^n,\psi_h)+(\textbf{g}^{n-\frac{1}{2}},\psi_h),
\end{alignat}
\end{subequations}
for all $(\textbf{v}_h,q_h,\textbf{H}_h,\psi_h)\in {\cV}^2_h\times \cM_h\times \cW^2_h\times \cQ^2_h$.
%
% Considering stabilized finite element pair, we revise the algorithm 5.1 as follows:
%
%\textbf{Algorithm 5.2:}  Given $(\textbf{u}^0_h,\textbf{B}^0_h,\textbf{w}^0_h) \in \cV^2_h\times \cW^2_h\times \cQ^2_h$,
%find $(\textbf{u}^n_h, p^n_h, \textbf{B}^n_h, \textbf{w}^n_h) \in {\cV}^2_h\times \cM^2_h\times \cW^2_h\times \cQ^2_h$, such that
% \begin{subequations}\label{eq61.1}
%\begin{alignat}{2} \label{eq61.1a}
% (\mathcal{D}\textbf{u}^{n+1}_h,\textbf{v}_h)+\overline{\mathbb{B}}(\overline{\textbf{u}}^n_h,\overline{p}^n_h;\textbf{v}_h,q_h)
%+b(\mathcal{I}(\textbf{u}^{n}_h),\overline{\textbf{u}}^n_h,\textbf{v}_h)\\\nonumber\quad
%+Sc_{\widehat{B}}(\mathcal{I}(\textbf{B}^{n}_h),\overline{\textbf{B}}^n_h,\textbf{v}_h)
%&=e(\overline{\textbf{w}}_h^n,\textbf{v})+(\textbf{f}^{n+\frac{1}{2}},\textbf{v}_h), \\
%(\mathcal{D}\textbf{B}^{n+1}_h,\textbf{H}_h)+a_B(\overline{\textbf{B}}^{n}_h,\textbf{H}_h)
%-c_{\widetilde{B}}(\overline{\textbf{u}}^{n}_h,\mathcal{I}(\textbf{B}^{n}_h),\textbf{H}_h)&=0, \\
%(\mathcal{D}\textbf{w}^{n+1}_h,\psi_h)+a_w(\overline{\textbf{w}}^{n}_h,\psi_h)+c_w(\mathcal{I}(\textbf{u}^{n}_h),\overline{\textbf{w}}^{n}_h,\psi_h)
%\\\nonumber\quad+2\nu_r(\overline{\textbf{w}}_h^n,\psi)&=e(\overline{\textbf{u}}_h^n,\psi)+(\textbf{g}^{n+\frac{1}{2}},\psi_h),
%\end{alignat}
%\end{subequations}
%for all $(\textbf{v}_h,q_h,\textbf{H}_h,\psi_h)\in  {\cV}^2_h\times \cM^2_h\times \cW^2_h\times \cQ^2_h$.

We may now state some unconditionally energy stable and optimal error estimates results for the second-order numerical scheme \eqref{eq60.1}.
\begin{theorem}  %Assume \textbf{Hypotheses A1-A3} hold.
Let $(\textbf{f},\textbf{g})\in L^2(0,T; H^{-1}(\Omega)^3)\times L^2(0,T; H^{-1}(\Omega)^3)$.
Then the finite element solution $(\textbf{u}^{n}_h,\textbf{B}^{n}_h,\textbf{w}^{n}_h)$ of problem \eqref{eq60.1} satisfy the following bounds
   \begin{align} \label{eq62} &\|\textbf{u}^N_h\|^2+S\|\textbf{B}^N_h\|^2
+\|\textbf{w}^N_h\|^2
\\\nonumber&\quad+{(\nu+\nu_r)}\Delta t\sum^{N}_{i=2}\|\nabla\overline{{\textbf{u}}}^n_h\|^2+c_{\star}\mu\Delta tS\sum^{N}_{i=2}\|\nabla\overline{{\textbf{B}}}^n_h\|^2\\\nonumber&\quad
+{\Delta t(c_a+c_d)}\sum^{N}_{i=2}\|\nabla\overline{{\textbf{w}}}^n_h\|^2+(c_0+c_d-c_a)\Delta t\sum^{N}_{i=2}\|\nabla\cdot\overline{{\textbf{w}}}^n_h\|^2\\\nonumber&
\leq C.
%\frac{\Delta t}{(\nu+\nu_r)}\sum^{N}_{i=2}\|\textbf{f}^{n-\frac{1}{2}}\|^2_{-1}+ \frac{\Delta t}{(c_a+c_d)}\sum^{N}_{i=2}\|\textbf{g}^{n-\frac{1}{2}}\|^2_{-1}+2\bigl(\|\textbf{u}^{0}_h\|^2+\|\textbf{u}^{1}_h\|^2\bigl)\\\nonumber&\quad
%+2S\bigl(\|\textbf{B}^{0}_h\|^2+\|\textbf{B}^{1}_h\|^2\bigr)
%+2\bigl(\|\textbf{w}^{0}_h\|^2+\|\textbf{w}^{1}_h\|^2\bigr).
 \end{align}
\end{theorem}

 \begin{theorem}   Suppose that \textbf{Assumpation A2} holds and the initial conditions satisfy
 \begin{align}\label{eq63}
\sum^1_{i=0}\bigl(\|\textbf{u}^i-\textbf{u}^i_h\|+ \|\textbf{B}^i-\textbf{B}^i_h\|+\|\textbf{w}^i-\textbf{w}^i_h\|\bigr)\leq Ch.
 \end{align}
Then the error estimates hold
\begin{align*}
\max_{2\leq n\leq N}\|\textbf{u}(t_n)-\textbf{u}^n_h\|^2
+\Delta t\sum_{n=2}^{N}\|\nabla (\textbf{u}(t_n)-\textbf{u}^n_h)\|^2 &\leq C\bigl(\Delta t^4+h^{2}\bigr),\\
\max_{2\leq n\leq N}\|\textbf{B}(t_n)-\textbf{B}^n_h\|^2
+\Delta t\sum_{n=2}^{N}\|\nabla (\textbf{B}(t_n)-\textbf{B}^n_h)\|^2 &\leq C\bigl(\Delta t^4+h^{2}\bigr),\\
\max_{2\leq n\leq N}\|\textbf{w}(t_n)-\textbf{w}^n_h\|^2
+\Delta t\sum_{n=2}^{N}\|\nabla (\textbf{w}(t_n)-\textbf{w}^n_h)\|^2 &\leq C\bigl(\Delta t^4+h^{2}\bigr).
 \end{align*}
\end{theorem}

 \begin{theorem}
 Suppose that \textbf{Assumpation A2} holds.  Suppose that the initial conditions satisfy \eqref{eq34} and
 \begin{subequations}\label{eq64.1}
\begin{alignat}{2} \label{eq64.1a}
\sum^1_{i=0}\bigl(\|\nabla(\textbf{u}^i-\textbf{u}^i_h)\|+\|\nabla(\textbf{B}^i-\textbf{B}^i_h)\|+\|\nabla(\textbf{w}^i-\textbf{w}^i_h)\|\bigr)\leq Ch,\\\label{eq42.1b}
d(\overline{\textbf{u}}^1_h,q)=0,\ \forall\, q\,\in \cM_h.
\end{alignat}
\end{subequations}
Then we have the following error estimates
\begin{align*}
\max_{2\leq n\leq N}\|\nabla(\textbf{u}(t_n)-\textbf{u}^n_h)\|^2
+\Delta t\sum_{n=2}^{N}\|\partial_t\textbf{u}(t_{n})-\mathcal{D}\textbf{u}^{n}_h\|^2 &\leq C\bigl(\Delta t^4+h^{2}\bigr),\\
\max_{2\leq n\leq N}\|\nabla(\textbf{B}(t_n)-\textbf{B}^n_h)\|^2
+\Delta t\sum_{n=2}^{N}\|\partial_t\textbf{B}(t_{n})-\mathcal{D}\textbf{B}^{n}_h\|^2 &\leq C\bigl(\Delta t^4+h^{2}\bigr),\\
\max_{2\leq n\leq N}\|\nabla(\textbf{w}(t_n)-\textbf{w}^n_h)\|^2
+\Delta t\sum_{n=2}^{N}\|\partial_t\textbf{w}(t_{n})-\mathcal{D}\textbf{w}^{n}_h\|^2 &\leq C\bigl(\Delta t^4+h^{2}\bigr).
 \end{align*}
\end{theorem}

\begin{theorem} Suppose that assumptions of Theorems 5.2-5.3 hold.
Then we have following error estimate
\begin{align*}
\Delta t\sum_{n=2}^{N}\|p(t_n)-p^n_h\|^2 &\leq C\bigl(\Delta t^4+h^{2}\bigr).
 \end{align*}
\end{theorem}

\textbf{Remark 3:}
%{\color{blue}}
Applying similar lines as the proof of {Theorem 3.1}, {Theorem 4.1}-{Theorem 4.3} and \cite{2019Ravindran},
 we can prove the results of {Theorem 5.1}-{Theorem 5.4}, here we omit it.

 \subsection{\label{Sec3} Decoupled numerical schemes }

In this section, we propose some fully discrete decoupled numerical schemes of problem \eqref{eq9.1} as follows.

\textbf{Scheme 5.3:}

 \emph{Step I:} Given $(\textbf{u}^{n-1}_h, \textbf{B}^{n-1}_h)\in \cV^1_h\times \cW^1_h$, for any $\textbf{H}_h\in \cW^1_h$, find $\textbf{B}^n_h\in \cW^1_h$, such that
\begin{align}\label{eq65}
(\mathcal{D}\textbf{B}^{n}_h,\textbf{H}_h)+a_B({\textbf{B}}^{n}_h,\textbf{H}_h)
-c_{\widetilde{B}}({\textbf{u}}^{n-1}_h,\textbf{B}^{n}_h,\textbf{H}_h)&=0.
 \end{align}

\emph{Step II:} Given $(\textbf{u}^{n-1}_h, \textbf{w}^{n-1}_h) \in \cV^1_h\times \cQ^1_h$, for any $\psi_h\in \cQ^1_h$, find $\textbf{w}^n_h\in \cW^1_h$, such that
 \begin{align}\label{eq66}
(\mathcal{D}\textbf{w}^n_h,\psi_h)+a_w(\textbf{w}^{n}_h,\psi_h)+c_w(\textbf{u}^{n-1}_h,\textbf{w}^{n}_h,\psi_h)
\\\nonumber\quad+2\nu_r(\textbf{w}_h^n,\psi_h)&=e(\textbf{u}_h^{n-1},\psi_h)+(\textbf{g}^n,\psi_h).
\end{align}

\emph{Step III:} Given $(\textbf{u}^{n-1}_h,\textbf{w}^{n-1}_h,\textbf{B}^{n-1}_h,\textbf{B}^n_h) \in \cV^1_h\times \cQ^1_h\times \cW^1_h\times \cW^1_h$,
for any $(\textbf{v}_h,q_h)\in \cV^1_h\times \cM_h$, find $(\textbf{u}^n_h,p^n_h)\in \cV^1_h\times \cM_h$,
such that
 \begin{subequations}\label{eq67.1}
\begin{alignat}{2} \label{eq67.1a}
 (\mathcal{D}\textbf{u}^n_h,\textbf{v}_h)+a_f(\textbf{u}^n_h,\textbf{v}_h)
+b(\textbf{u}^{n-1}_h,\textbf{u}^n_h,\textbf{v}_h)\\\nonumber\quad
+Sc_{\widehat{B}}(\textbf{B}^{n-1}_h,\textbf{B}^n_h,\textbf{v}_h)-d(\textbf{v}_h,p^n_h)
&=e(\textbf{w}_h^{n-1},\textbf{v}_h)+(\textbf{f}^n,\textbf{v}_h), \\
d(\textbf{u}^n_h,q_h)&=0.
\end{alignat}
\end{subequations}

\textbf{Scheme 5.4:}

 \emph{Step I:} Given $(\textbf{u}^{n-1}_h, \textbf{B}^{n-1}_h)\in \cV^2_h\times \cW^2_h$, for any $\textbf{H}_h\in \cW^2_h$, find $\textbf{B}^n_h\in \cW^2_h$, such that
 \begin{align}\label{eq68}
(\mathcal{D}\textbf{B}^{n}_h,\textbf{H}_h)+a_B({\textbf{B}}^{n}_h,\textbf{H}_h)
-c_{\widetilde{B}}({\textbf{u}}^{n-1}_h,\textbf{B}^{n}_h,\textbf{H}_h)&=0.
 \end{align}

\emph{Step II:} Given $(\textbf{u}^{n-1}_h, \textbf{w}^{n-1}_h) \in \cV^2_h\times \cQ^2_h$, for any $\psi_h\in \cQ^2_h$, find $\textbf{w}^n_h\in \cW^2_h$, such that
 \begin{align}\label{eq69}
(\mathcal{D}\textbf{w}^n_h,\psi_h)+a_w(\textbf{w}^{n}_h,\psi_h)+c_w(\textbf{u}^{n-1}_h,\textbf{w}^{n}_h,\psi_h)
\\\nonumber\quad+2\nu_r(\textbf{w}_h^n,\psi_h)&=e(\textbf{u}_h^{n-1},\psi_h)+(\textbf{g}^n,\psi_h).
\end{align}

\emph{Step III:} Given $(\textbf{u}^{n-1}_h,\textbf{w}^{n-1}_h,\textbf{B}^{n-1}_h,\textbf{B}^n_h) \in \cV^2_h\times \cQ^2_h\times \cW^2_h\times \cW^2_h$,
for any $(\textbf{v}_h,q_h)\in \cV^2_h\times \cM_h$, find $(\textbf{u}^n_h,p^n_h)\in \cV^2_h\times \cM_h$,
such that
 \begin{subequations}\label{eq70.1}
\begin{alignat}{2} \label{eq70.1a}
 (\mathcal{D}\textbf{u}^n_h,\textbf{v}_h)+\overline{\mathbb{B}}(\textbf{u}^n_h,p^n_h;\textbf{v}_h,q_h)
+b(\textbf{u}^{n-1}_h,\textbf{u}^n_h,\textbf{v}_h)\\\nonumber\quad
+Sc_{\widehat{B}}(\textbf{B}^{n-1}_h,\textbf{B}^n_h,\textbf{v}_h)
&=e(\textbf{w}_h^{n-1},\textbf{v}_h)+(\textbf{f}^n,\textbf{v}_h).
\end{alignat}
\end{subequations}

It can be noted that the calculation of all nonlinear terms will not bring any variable coefficients. In each
time step, the total computational cost of the decoupled schemes is to solve several independent, general maxwell equations,
elliptic equations and stokes
equations.
Because the nonlinear terms $c_{\widetilde{B}}({\textbf{u}}^{n-1}_h,\textbf{B}^{n}_h,\textbf{H}_h)$
and $c_{\widehat{B}}(\textbf{B}^{n-1}_h,\textbf{B}^n_h,\textbf{v}_h)$
are unable to be controlled,  we can't obtain similar stability results for the Schemes $5.3$ and $5.4$.
 We will implement
the stability and error analysis in the future work by following the similar lines  in \cite{2022Chen,He2015,2017Hasler,2022Wang}.

\textbf{Remark 4:}
We note that the pressure of the Schemes 3.1, 5.1 and 5.3 is implicitly updated.
If a decoupled Stokes solver is introduced, using the similar lines in  \cite{2015Chen,2019Chen,2022Chen,2022Wang}
the energy stability and convergence analysis can be proved.

\section{\label{Sec5} Numerical results }

In this section, some numerical tests are given to confirm
 the rates of convergence and energy stable
  of the first/second-order numerical Schemes 3.1-3.2 and 5.1-5.4.
 The tests have been solved by
using the software package \emph{FreeFem++ \cite{Hecht2012}}.

 The conforming finite element pair $\textbf{P}_{1b}-P_{1}-\textbf{P}_{1b}-\textbf{P}_{1b}$ and
 stabilized finite element pair $\textbf{P}_{1}-P_{1}-\textbf{P}_{1}-\textbf{P}_{1}$ for the velocity, the pressure, the magnetic and the micro-rotation are employed.
We consider the computing domain to be $\Omega=[0,1]^3$.  The uniform triangles and tetrahedra meshes are utilized.
For simplicity, we denote the following discrete errors:
 \begin{align*}
  E(\textbf{u})&=:\max_{1\leq n\leq N}\|\nabla(\textbf{u}(t_n)-\textbf{u}^n_h)\|,\qquad &&E(\textbf{B})=:\max_{1\leq n\leq N}\|\nabla(\textbf{B}(t_n)-\textbf{B}^n_h)\|,\\
  E(\textbf{w})&=:\max_{1\leq n\leq N}\|\nabla(\textbf{w}(t_n)-\textbf{w}^n_h)\|,\qquad &&
  E(p)=:\Bigl(\Delta t\sum^N_{n=1}\|p(t_n)-p^n_h\|^2\Bigr)^{\frac{1}{2}}.
 \end{align*}

\subsection{\label{Sec5} 3D convergence test}
The known source terms, initial conditions and boundary conditions are taken
so that the exact solution is defined by
\begin{equation*}
\begin{array}{ll}
u_1(x,y,z,t)&= -2(1-\cos(2\pi x))\sin(2\pi y)\sin(2\pi z)\cos(t),\\
u_2(x,y,z,t)&=(1-\cos(2\pi y))\sin(2\pi x)\sin(2\pi z)\cos(t),\\
u_3(x,y,z,t)&=(1-\cos(2\pi z))\sin(2\pi x)\sin(2\pi y)\cos(t),\\
p(x,y,z,t)&= (\sin(4\pi x)+\sin(4\pi y)+\sin(4\pi z))\cos(t),\\
B_1(x,y,z,t)&=  2\sin(\pi x)\cos(\pi y)\cos(\pi z)\cos(t),\\
B_2(x,y,z,t)&=-\sin(\pi y)\cos(\pi x)\cos(\pi z)\cos(t),\\
B_3(x,y,z,t)&=-\sin(\pi z)\cos(\pi x)\cos(\pi y)\cos(t),\\
w_1(x,y,z,t)&= (1-\cos(2\pi x))\sin(2\pi y)\sin(2\pi z)\cos(t),\\
w_2(x,y,z,t)&=(1-\cos(2\pi y))\sin(2\pi x)\sin(2\pi z)\cos(t),\\
w_3(x,y,z,t)&=(1-\cos(2\pi z))\sin(2\pi x)\sin(2\pi y)\cos(t).
\end{array}\end{equation*}
%as the  exact solution, and some known source terms are chosen such that
%the  exact solution satisfy the problem \eqref{eq1.1}-\eqref{eq3}.

 The performance of the coupled/decoupled numerical Schemes 3.1-3.2 and 5.1-5.4
 herein are studied to present the accuracy
experiments with respect to $\Delta t$ and $h$.
The parameters are taken to $\nu_r=1$, $\nu=1$, $\mu=1$, $c_0= c_a=c_d=0.5$ and $S= 1$.
The terminal time is chosen to $T=1$. The time steps are set to $\Delta t=h$ when Schemes 3.1-3.2
and 5.3-5.4, and $\Delta t=\sqrt{h}$ when Schemes 5.1-5.2.
We take the spatial grid size $h=\frac{1}{4},\frac{1}{8},\frac{1}{12},\frac{1}{16}$.
 From tables 1-4, we can see that the rates of error estimate for Schemes 3.1-3.2 and 5.1-5.2 are consistent
with the theoretical results. In addition, the error
rates of the velocity, the magnetic and the micro-rotation can achieve optimal.
 We get an interesting observation that the error of the pressure has
better convergence rate than theoretical results.
From tables 5-6, it is easy to present that the rates of convergence in
time and space are all first order for the Schemes 5.3-5.4.
\begin{center}
\begin{tabular}{ccccccccccccc}
\hline ${\small h}$& ${\small  E(\textbf{u})}$ &${\small  Rate }$& ${\small E(\textbf{B}) }$
 &   ${\small  Rate }$ & ${\small E(\textbf{w})}$ &${\small  Rate }$ & ${\small E(p)}$ &${\small  Rate }$ \\
\hline
$\frac{1}{4}$& 10.4061& /&     2.08861&/  &7.04381 &/ & 17.4691& /\\
$\frac{1}{8}$&   6.05154  &  0.7821&1.14118 & 0.8720& 4.20652 & 0.7437&  6.85828&1.3489\\
$\frac{1}{12}$&   4.14405&0.9338& 0.774757  &  0.9551&  2.90427& 0.9137&3.6057 & 1.5857\\
$\frac{1}{16}$& 3.13666 & 0.9681&   0.584914 & 0.9771& 2.20505&0.9574& 2.29992& 1.5630\\
\hline
\end{tabular}
\end{center}
\centerline{Table 1\quad Numerical results of the with $\Delta t=h$  and $\textbf{P}_{1b}-P_1-\textbf{P}_{1b}-\textbf{P}_{1b}$ for Scheme 3.1.}
\begin{center}
\begin{tabular}{ccccccccccccc}
\hline ${\small h}$& ${\small  E(\textbf{u})}$ &${\small  Rate }$& ${\small E(\textbf{B}) }$
 &   ${\small  Rate }$ & ${\small E(\textbf{w})}$ &${\small  Rate }$ & ${\small E(p)}$ &${\small  Rate }$ \\
\hline$\frac{1}{4}$& 10.0961& /&      2.11799 &/  &7.15628   &/&3.40901 &/\\
$\frac{1}{8}$& 6.01581 &  0.7470&1.15554&0.8741& 4.2608  &0.7481&1.57102 & 1.1177\\
$\frac{1}{12}$& 4.15617 & 0.9120&   0.784304& 0.9558& 2.94003&  0.9151&   0.838178& 1.5495\\
$\frac{1}{16}$& 3.15621 &   0.9567& 0.592073&   0.9774&  2.23176 &  0.9581& 0.513859 &1.7008\\
\hline
\end{tabular}
\end{center}
\centerline{Table 2\quad Numerical results of the with $\Delta t=h$  and $\textbf{P}_{1}-P_1-\textbf{P}_{1}-\textbf{P}_{1}$ for Scheme 3.2.}
\begin{center}
\begin{tabular}{ccccccccccccc}
\hline ${\small h}$& ${\small  E(\textbf{u})}$ &${\small  Rate }$& ${\small E(\textbf{B}) }$
 &   ${\small  Rate }$ & ${\small E(\textbf{w})}$ &${\small  Rate }$& ${\small E(p)}$ &${\small  Rate }$ \\
\hline$\frac{1}{4}$& 15.2988& /&     3.27729&/  &10.2113&/ & 25.4014&/\\
$\frac{1}{8}$&   10.3227& 0.5676 & 1.97734 &0.7289&7.07013&0.5304&9.55111 &1.4112\\
$\frac{1}{12}$& 7.27653&0.8624& 1.36738&0.9097& 5.02573&0.8418 &4.46046 & 1.8779\\
$\frac{1}{16}$& 5.62094& 0.8974&  1.04924& 0.9206& 3.89457& 0.8864& 3.00632&1.3714\\
\hline
\end{tabular}
\end{center}
\centerline{Table 3\quad Numerical results of the with $\Delta t=\sqrt{h}$  and $\textbf{P}_{1b}-P_1-\textbf{P}_{1b}-\textbf{P}_{1b}$ for Scheme 5.1.}
\begin{center}
\begin{tabular}{ccccccccccccc}
\hline ${\small h}$& ${\small  E(\textbf{u})}$ &${\small  Rate }$& ${\small E(\textbf{B}) }$
 &   ${\small  Rate }$ & ${\small E(\textbf{w})}$ &${\small  Rate }$ & ${\small E(p)}$ &${\small  Rate }$ \\
\hline$\frac{1}{4}$&6.4329& /&    1.30688&/  &4.30119 &/& 2.76459&/\\
$\frac{1}{8}$& 4.73988 &0.4406&0.909978&0.5222& 3.28828&0.3874&1.41325 &0.9680\\
$\frac{1}{12}$&3.53837 &0.7210&  0.667768& 0.7633&  2.4774& 0.6983&0.771974 &1.4914 \\
$\frac{1}{16}$&2.79563& 0.8190&0.52423&0.8412&  1.9641&0.8071 &  0.525558& 1.3365\\
\hline
\end{tabular}
\end{center}
\centerline{Table 4\quad Numerical results of the with $\Delta t=\sqrt{h}$  and $\textbf{P}_{1}-P_1-\textbf{P}_{1}-\textbf{P}_{1}$ for Scheme 5.2.}
\begin{center}
\begin{tabular}{ccccccccccccc}
\hline ${\small h}$& ${\small  E(\textbf{u})}$ &${\small  Rate }$& ${\small E(\textbf{B}) }$
 &   ${\small  Rate }$ & ${\small E(\textbf{w})}$ &${\small  Rate }$ & ${\small E(p)}$ &${\small  Rate }$ \\
\hline$\frac{1}{4}$& 10.406& /&   2.08865&/  & 7.03625&/&  17.4771&/\\
$\frac{1}{8}$& 6.05152&0.7820&1.14077 &0.8726&  4.20326& 0.7433&6.85977&1.3492\\
$\frac{1}{12}$&  4.14405&0.9338&  0.774602&0.9547&2.90306&0.9128& 3.60628&1.5858\\
$\frac{1}{16}$&3.13667& 0.9681& 0.584848 & 0.9768& 2.20452&0.9568&2.30023 &1.5631\\
\hline
\end{tabular}
\end{center}
\centerline{Table 5\quad Numerical results of the with $\Delta t=h$  and $\textbf{P}_{1b}-P_1-\textbf{P}_{1b}-\textbf{P}_{1b}$ for Scheme  5.3.}
\begin{center}
\begin{tabular}{ccccccccccccc}
\hline ${\small h}$& ${\small  E(\textbf{u})}$ &${\small  Rate }$& ${\small E(\textbf{B}) }$
 &   ${\small  Rate }$ & ${\small E(\textbf{w})}$ &${\small  Rate }$ & ${\small E(p)}$ &${\small  Rate }$ \\
\hline$\frac{1}{4}$&10.0959& /&   2.11813&/  & 7.15241&/&3.41106 &/\\
$\frac{1}{8}$&  6.01576& 0.7470& 1.15533 & 0.8745&  4.2586& 0.7481&1.57151 &1.1181\\
$\frac{1}{12}$&   4.15616&0.9120& 0.784219& 0.9556&2.93914&0.9146&0.838513 &1.5492\\
$\frac{1}{16}$&3.15622&0.9567&0.592037 &0.9772&2.23137 &0.9577& 0.514115& 1.7004\\
\hline
\end{tabular}
\end{center}
\centerline{Table 6\quad Numerical results of the with $\Delta t=h$  and $\textbf{P}_{1}-P_1-\textbf{P}_{1}-\textbf{P}_{1}$ for Scheme 5.4.}
%\begin{center}
%\begin{tabular}{ccccccccccccc}
%\hline ${\small h}$& ${\small  E(\textbf{u})}$ &${\small  Rate }$& ${\small E(\textbf{B}) }$
% &   ${\small  Rate }$ & ${\small E(\textbf{w})}$ &${\small  Rate }$ & ${\small E(p)}$ &${\small  Rate }$ \\
%\hline$\frac{1}{4}$&15.2954& /&    3.27726&/  & 10.199&/& 25.4249&/\\
%$\frac{1}{8}$& 10.321& 0.5675&1.97438&0.7311& 7.06308& 0.5301& 9.53877 & 1.4144\\
% $\frac{1}{12}$& 7.09323&0.9250&1.33221& 0.9703&4.89804&0.9028& 4.9037&1.6410\\
%$\frac{1}{16}$&5.62061& 0.8089& 1.04861&0.8321&  3.89328&  0.7980&3.0002 &1.7078\\
%\hline
%\end{tabular}
%\end{center}
%\centerline{Table 7\quad Numerical results of the with $\Delta t=\sqrt{h}$  and $\textbf{P}_{1b}-P_1-\textbf{P}_{1b}-\textbf{P}_{1b}$ for algorithm 6.3.}
%\begin{center}
%\begin{tabular}{ccccccccccccc}
%\hline ${\small h}$& ${\small  E(\textbf{u})}$ &${\small  Rate }$& ${\small E(\textbf{B}) }$
% &   ${\small  Rate }$ & ${\small E(\textbf{w})}$ &${\small  Rate }$ & ${\small E(p)}$ &${\small  Rate }$ \\
%\hline$\frac{1}{4}$&6.42913& /&  1.30719&/  &  4.28292&/&2.7598 &/\\
%$\frac{1}{8}$&  4.88998 &0.3948&0.935275&0.4830&3.38552&  0.3392& 1.40944&0.9694\\
%$\frac{1}{12}$& 3.32243  & 0.9532&  0.625359& 0.9927& 2.31866& 0.9335&0.793825 &1.4159\\
%$\frac{1}{16}$& 2.79517& 0.6007&0.523413& 0.6186& 1.96145&0.5816& 0.523082& 1.4500\\
%\hline
%\end{tabular}
%\end{center}
%\centerline{Table 8\quad Numerical results of the with $\Delta t=\sqrt{h}$  and $\textbf{P}_{1}-P_1-\textbf{P}_{1}-\textbf{P}_{1}$ for algorithm 6.4.}

\subsection{\label{Sec5} 3D stability test}

In this test, we present the unconditionally energy stable of our established Schemes 3.1-3.2 and 5.1-5.4 for $3$D magneto-micropolar fluid flows.
The dissipative property of  $3$D magneto-micropolar equations \eqref{eq1.1}-\eqref{eq3} is studied. To this goal, we take
$\textbf{f}=\textbf{g}=\textbf{0}$ and set the initial conditions as follows:
\begin{equation*}
\begin{array}{ll}
u_1(x,y,z,0)&= -2(1-\cos(2\pi x))\sin(2\pi y)\sin(2\pi z),\\
u_2(x,y,z,0)&=(1-\cos(2\pi y))\sin(2\pi x)\sin(2\pi z),\\
u_2(x,y,z,0)&=(1-\cos(2\pi z))\sin(2\pi x)\sin(2\pi y),\\
%p(x,y,z,0)&= (\sin(4\pi x)+\sin(4\pi y)+\sin(4\pi z)),\\
B_1(x,y,z,0)&=  2\sin(\pi x)\cos(\pi y)\cos(\pi z),\\
B_2(x,y,z,0)&=-\sin(\pi y)\cos(\pi x)\cos(\pi z),\\
B_3(x,y,z,0)&=-\sin(\pi z)\cos(\pi x)\cos(\pi y),\\
w_1(x,y,z,0)&= (1-\cos(2\pi x))\sin(2\pi y)\sin(2\pi z),\\
w_2(x,y,z,0)&=(1-\cos(2\pi y))\sin(2\pi x)\sin(2\pi z),\\
w_2(x,y,z,0)&=(1-\cos(2\pi z))\sin(2\pi x)\sin(2\pi y).
\end{array}\end{equation*}
%We consider computing domain to be $\Omega=[0,1]^3$,
 The parameters are also set to $\nu_r=1$, $\nu=0.001$, $\mu=0.01$, $c_0= c_a=0.005,c_d=0.5$ and $S= 1$.
In this test, we take the terminal $T=12$ and mesh $h=\frac{1}{13}$.
By using the discrete energy $E_n=\|\textbf{u}_h^n\|^2+\|\textbf{B}_h^n\|^2+\|\textbf{w}_h^n\|^2$,
we compare the long time unconditionally stable of our proposed Schemes 3.1-3.2, Schemes 5.1-5.4.
 In Figs. 1-3, we show the time evolution of
the discrete energy $E_n$ for time step sizes $\Delta t=0.125,0.25,0.5$ until $T = 12$. It is easy to show that
three curves present monotonic decay, which numerically verify that our schemes are
unconditionally energy stable.

\begin{figure}
\centering
\subfigure[]{
\includegraphics[scale=0.35]{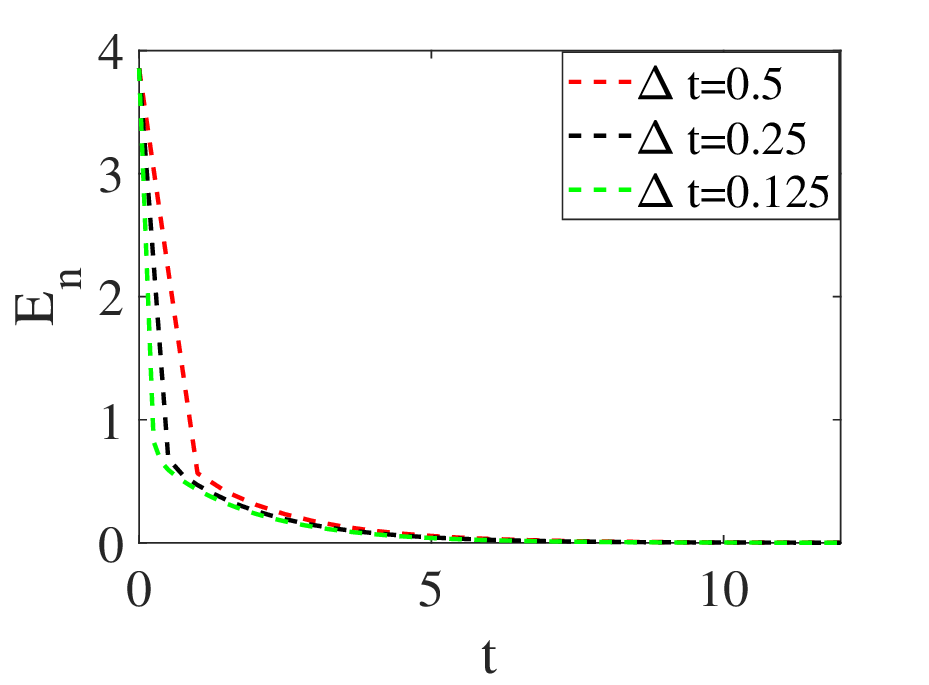}}
\subfigure[]{
\includegraphics[scale=0.35]{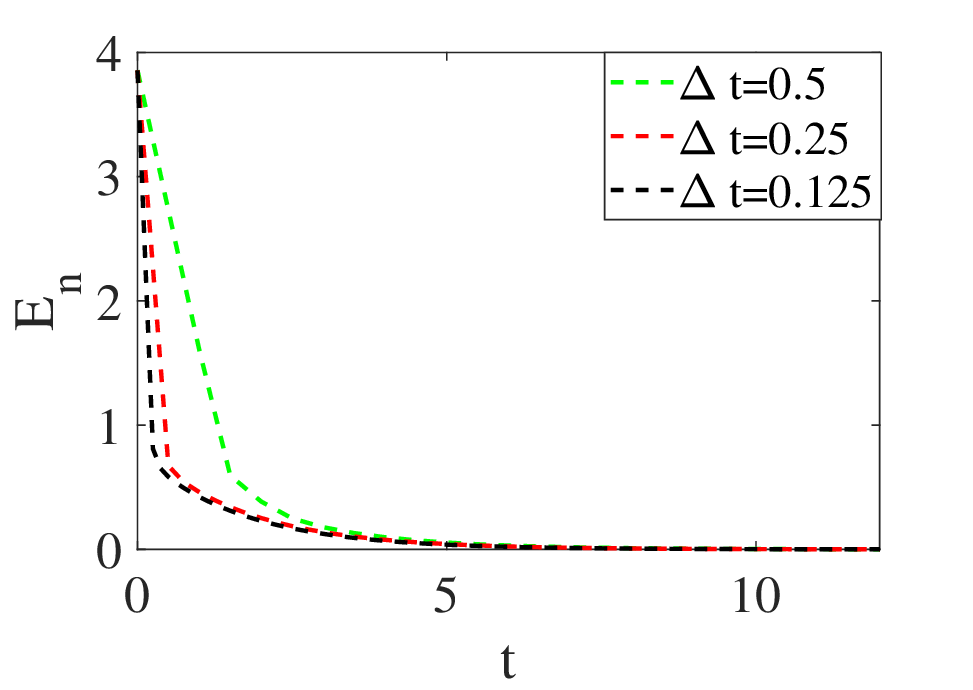}}
\caption{(a) Discrete energy $E_n$ with Scheme  3.1; (b) Discrete energy $E_n$ with Scheme  3.2.}
\label{fig1}
\end{figure}

\begin{figure}
\centering
\subfigure[]{
\includegraphics[scale=0.35]{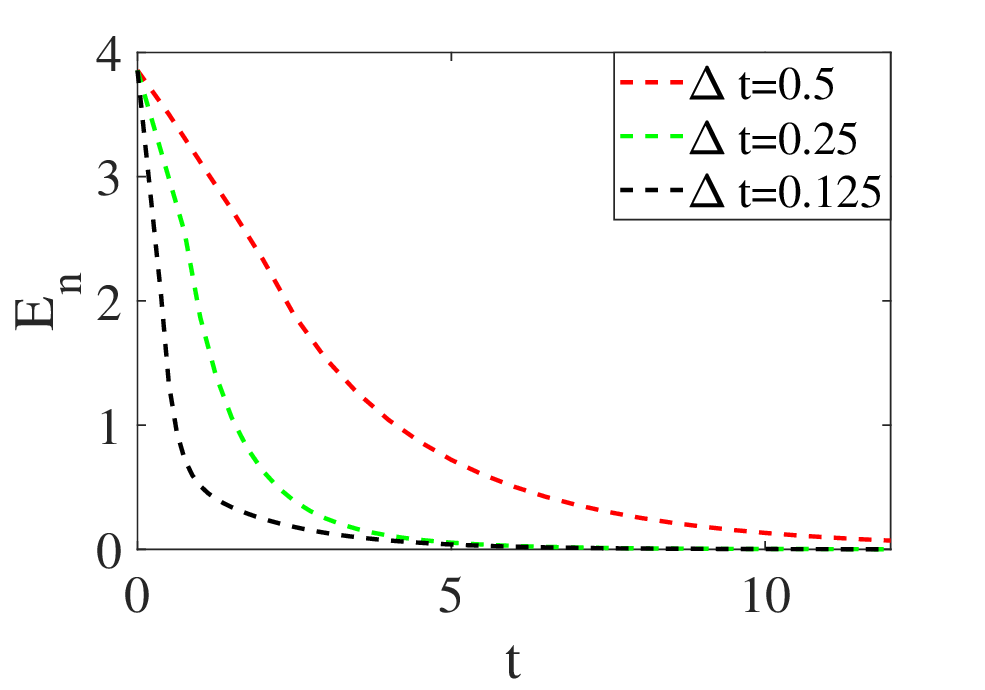}}
\subfigure[]{
\includegraphics[scale=0.35]{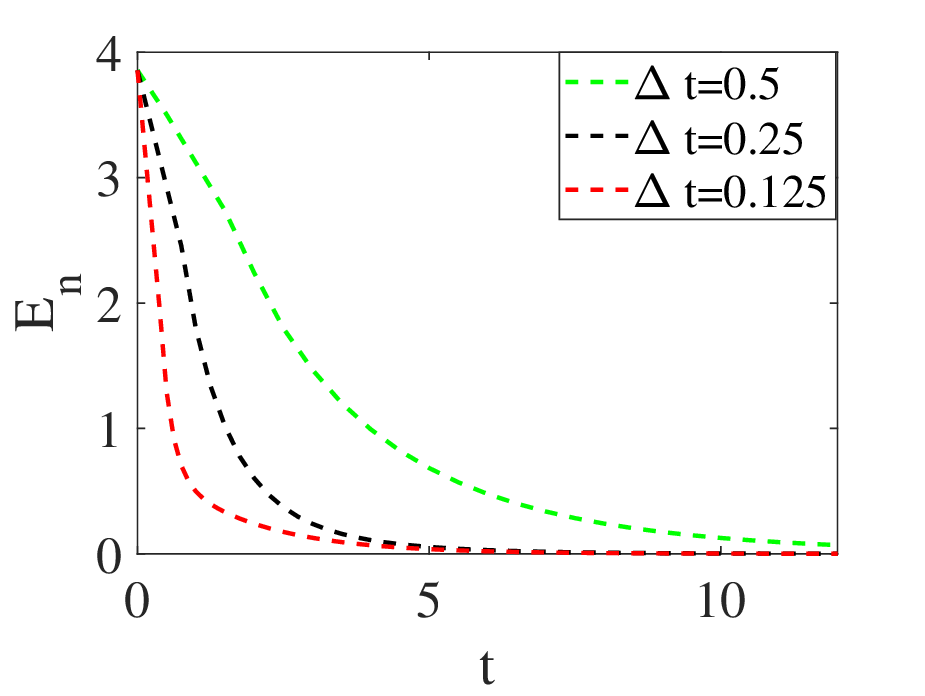}}
\caption{(a) Discrete energy $E_n$ with Scheme 5.1; (b) Discrete energy $E_n$ with Scheme 5.2.}
\label{fig1}
\end{figure}

\begin{figure}
\centering
\subfigure[]{
\includegraphics[scale=0.35]{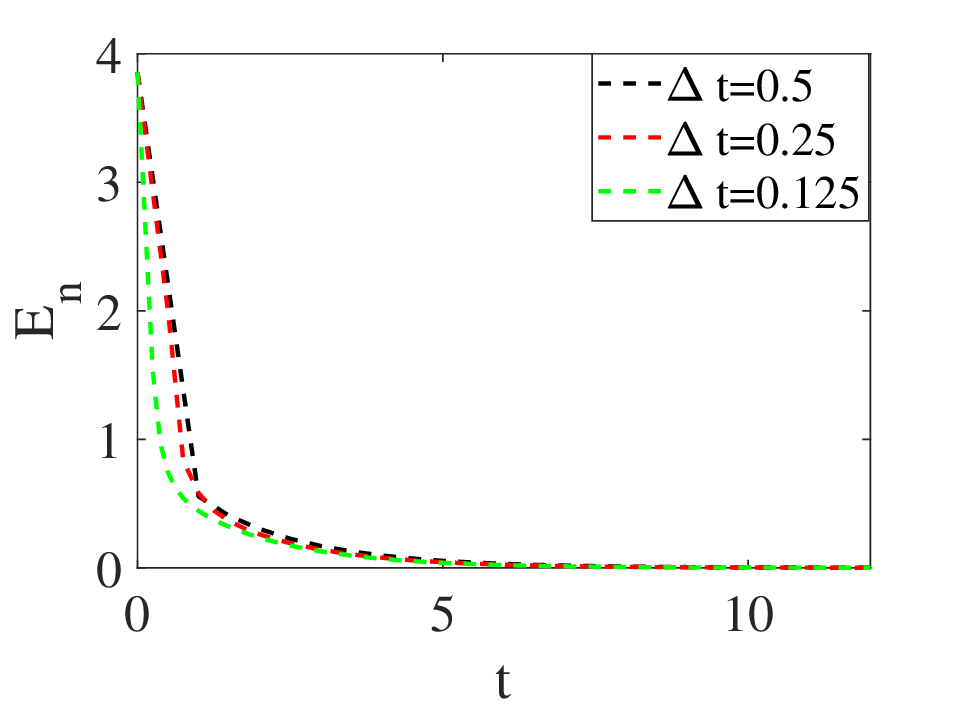}}
\subfigure[]{
\includegraphics[scale=0.35]{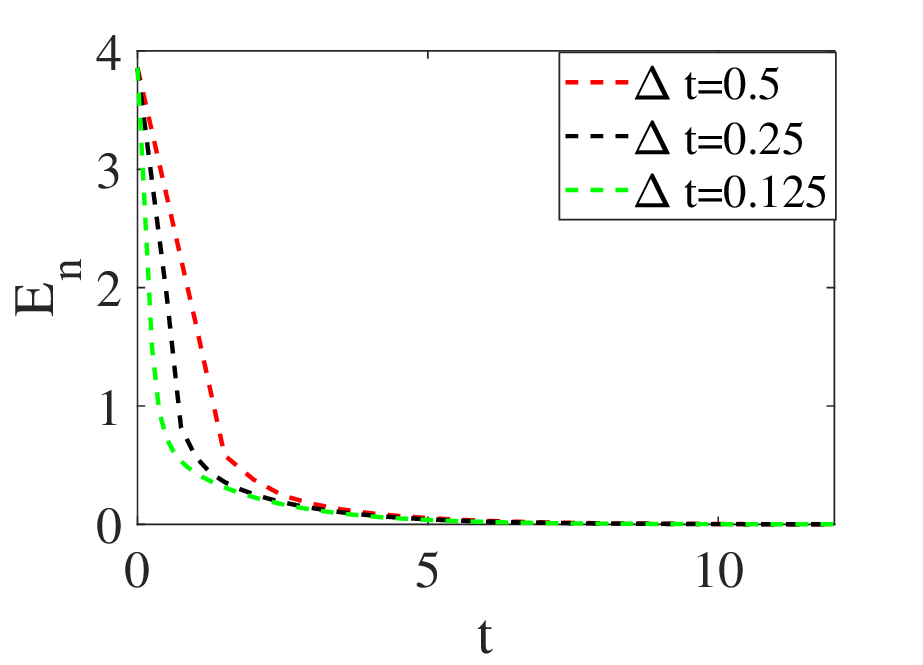}}
\caption{(a) Discrete energy $E_n$ with Scheme 5.3; (b) Discrete energy $E_n$ with Scheme 5.4.}
\label{fig1}
\end{figure}

%\begin{figure}
%\centering
%\subfigure[]{
%\includegraphics[scale=0.35]{3D_second_p1b_decouped.eps}}
%\subfigure[]{
%\includegraphics[scale=0.35]{3D_second_p1_decouped.eps}}
%\caption{(a) Discrete energy $E_n$ with algorithm 6.3; (b) Discrete energy $E_n$ with algorithm 6.4.}
%\label{fig1}
%\end{figure}

\subsection{\label{Sec5.3} 3D Lid-driven cavity flow }

In this test, we show the long time stability of 3D Lid-driven cavity flow simulation.
   The parameters are taken to $\nu_r=0.1$, $\nu=0.1$, $\mu=0.1$, $c_0= c_a=c_d=0.5$ and $S= 1$.
  % We consider this example on a unit volume $\Omega=[0,1]^3$.
 The boundary conditions are given as follows:
\begin{equation*}
\begin{array}{ll}
\textbf{u}=(0,0,0),&\mbox{on}\ \ x=0,1, \ \ y=0,1\ \mbox{and} \ z=0,\\
\textbf{u}=(1,0,0),&\mbox{on}\ \ z=1,\\
\textbf{n}\times \textbf{B}=\textbf{n}\times \textbf{B}_d, &\mbox{on}\ \  \partial\Omega,\\
\textbf{w}=(0,0,0),&\mbox{on}\ \ x=0,1, \ \ y=0,1\ \mbox{and} \ z=0,\\
\textbf{w}=(0,0,1),&\mbox{on}\ \ z=1,\\
\end{array}\end{equation*}
where $\textbf{B}_d = (0,0,1)$.

The mesh size and time step size are chosen to $\Delta t=0.05$ and $h=\frac{1}{12}$, respectively.
Numerical solutions were computed by utilizing Scheme 5.1 and algorithm 5.3 at $T = 5$ are shown in Figs.4-6.
We can see that our proposed Scheme 5.1 and Scheme 5.3 satisfy the long time unconditionally stable for 3D Lid-driven cavity flow.

\begin{figure}
\centering
\subfigure[]{
\includegraphics[scale=0.25]{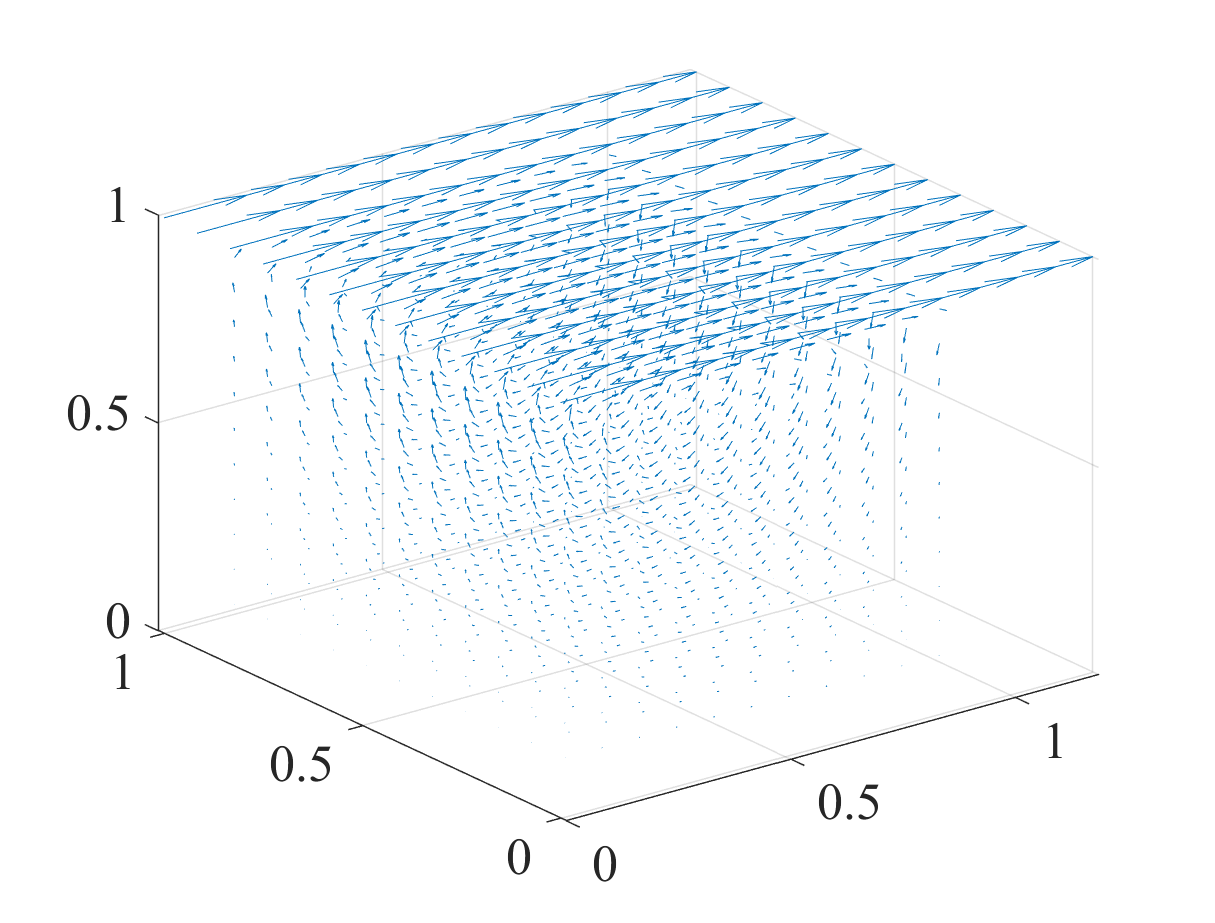}}
\subfigure[]{
\includegraphics[scale=0.25]{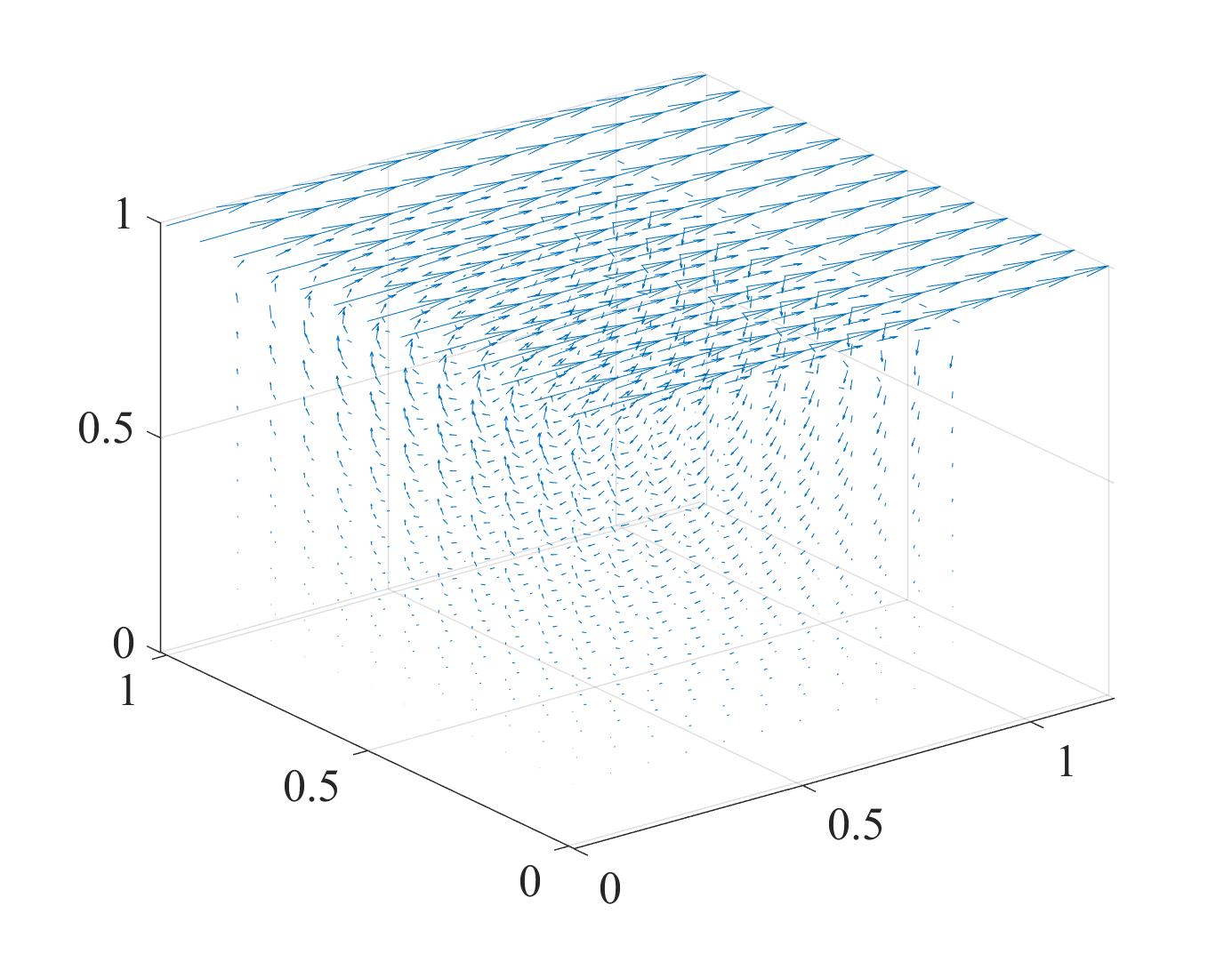}}
\caption{\label{fig:0}\small  The velocity field (a) and (b) of flow by the Scheme  5.1 and Scheme 5.3.}
\end{figure}

\begin{figure}
\centering
\subfigure[]{
\includegraphics[scale=0.23]{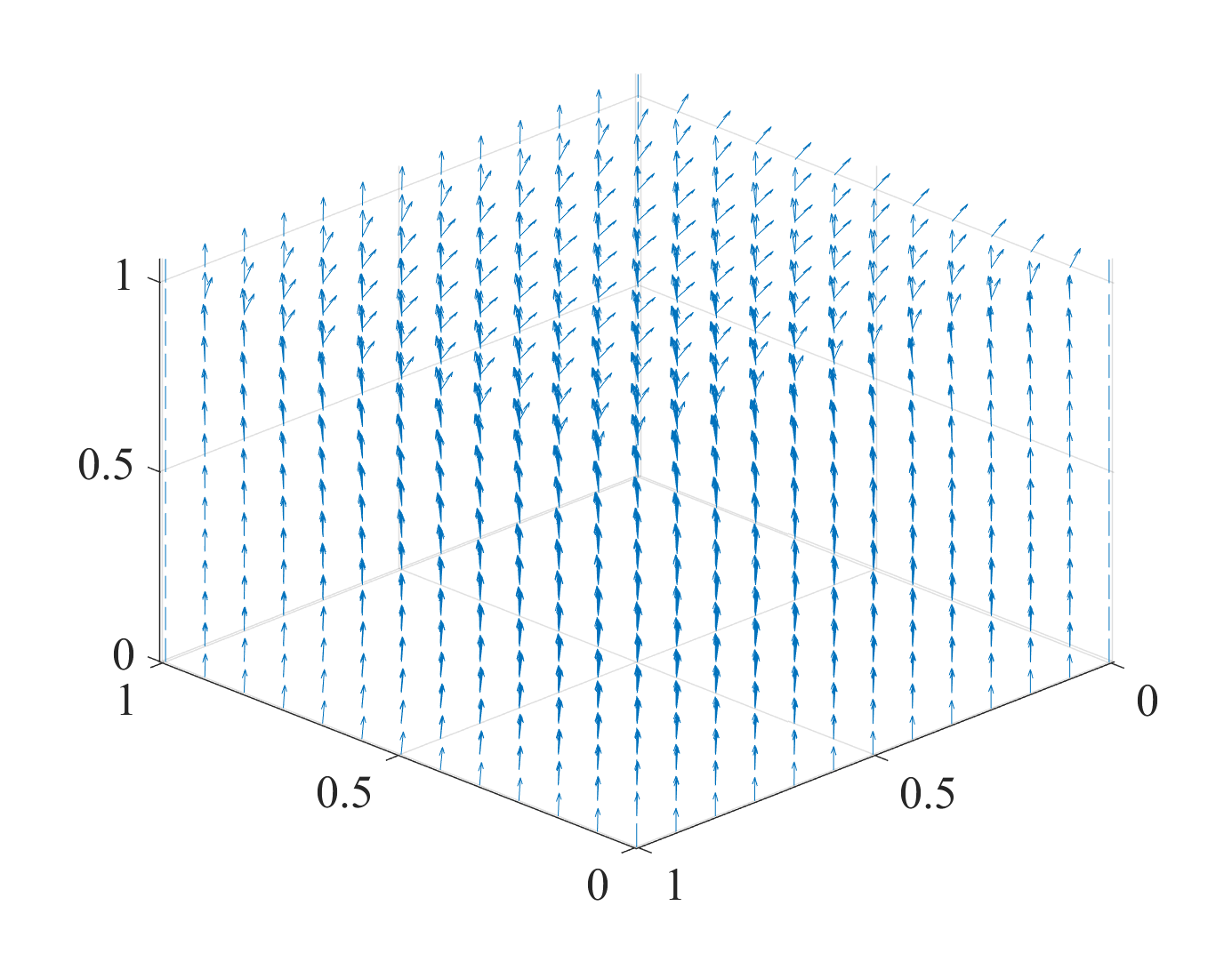}}
\subfigure[]{
\includegraphics[scale=0.23]{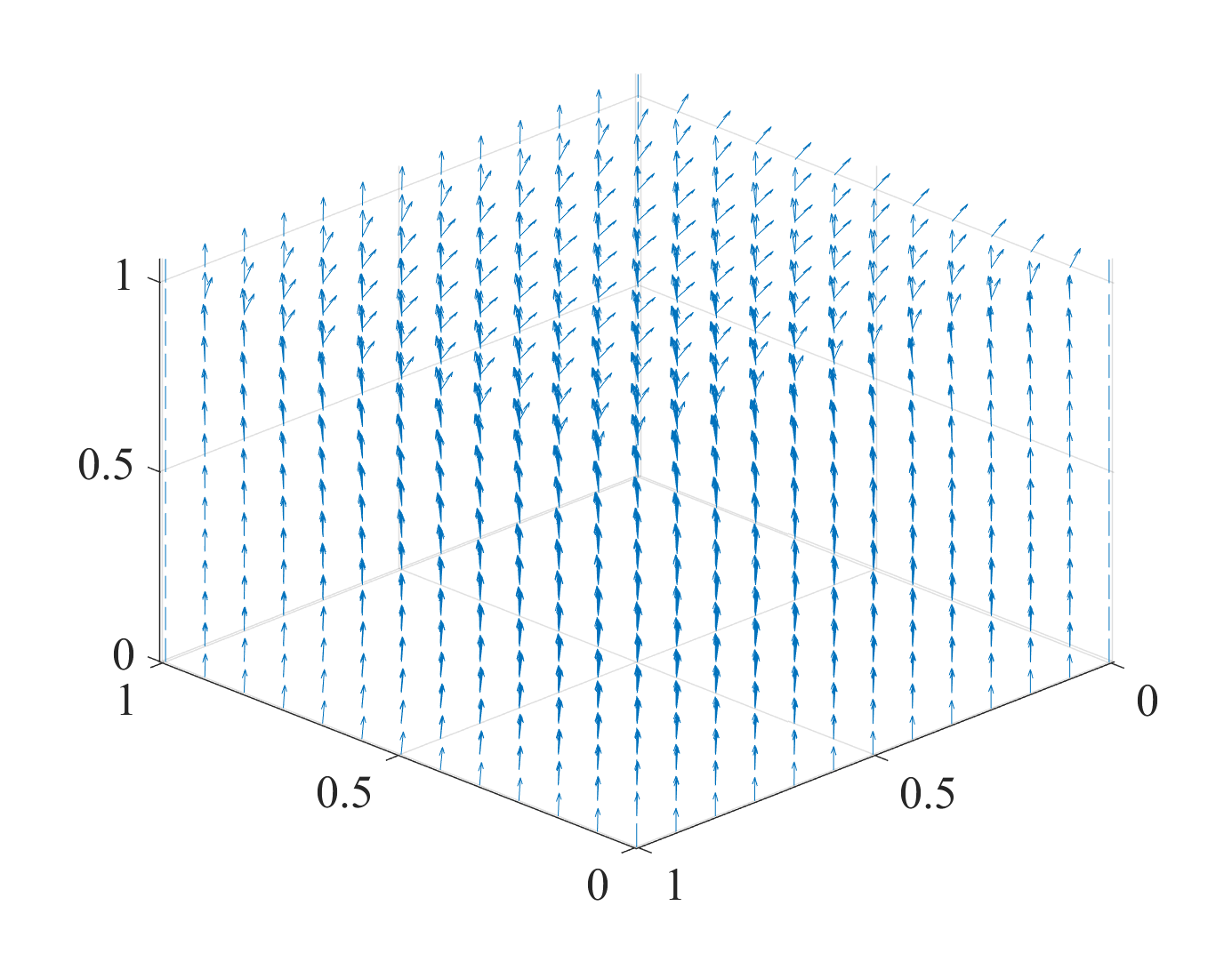}}
\caption{\label{fig:0}\small  The magnetic field (a) and (b) of flow by the Scheme 5.1 and Scheme 5.3.}
\end{figure}

\begin{figure}
\centering
\subfigure[]{
\includegraphics[scale=0.25]{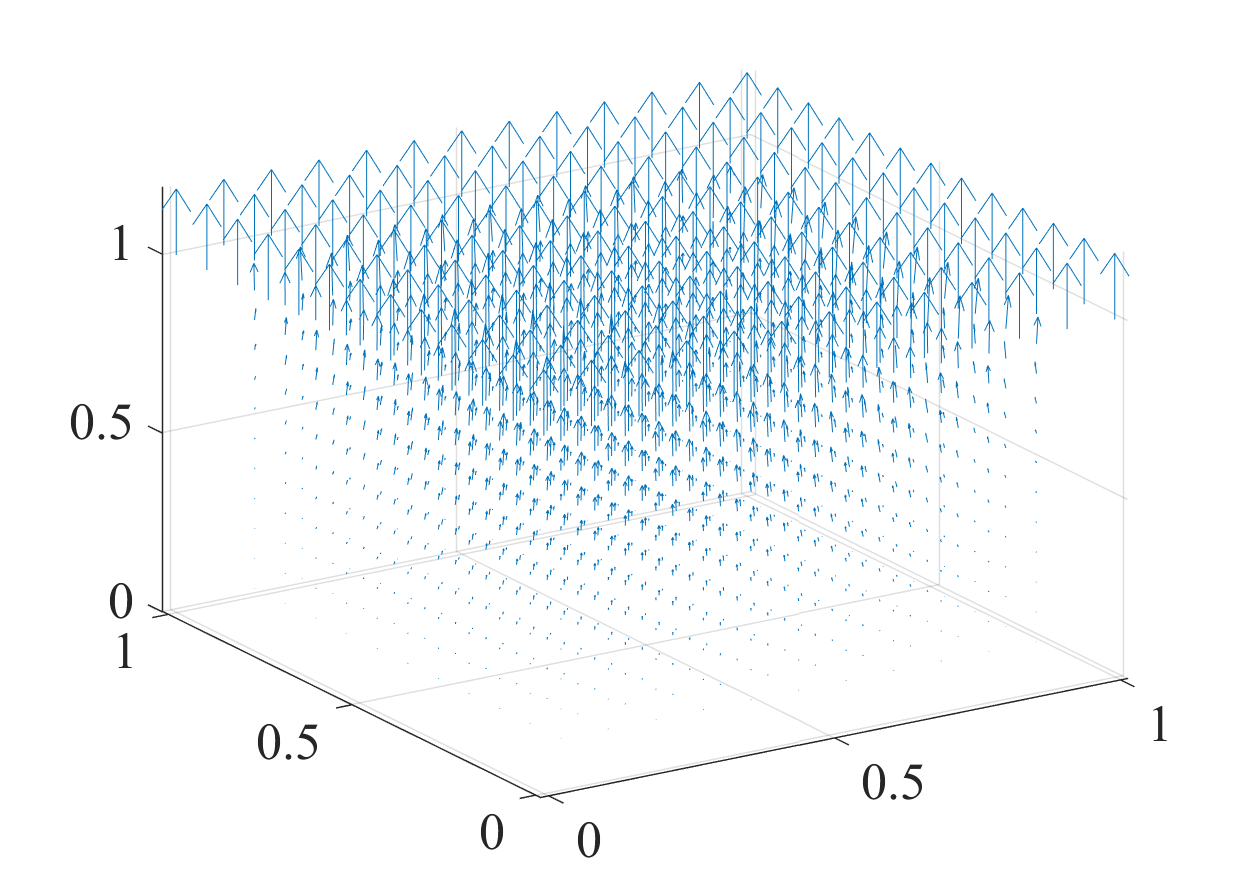}}
\subfigure[]{
\includegraphics[scale=0.25]{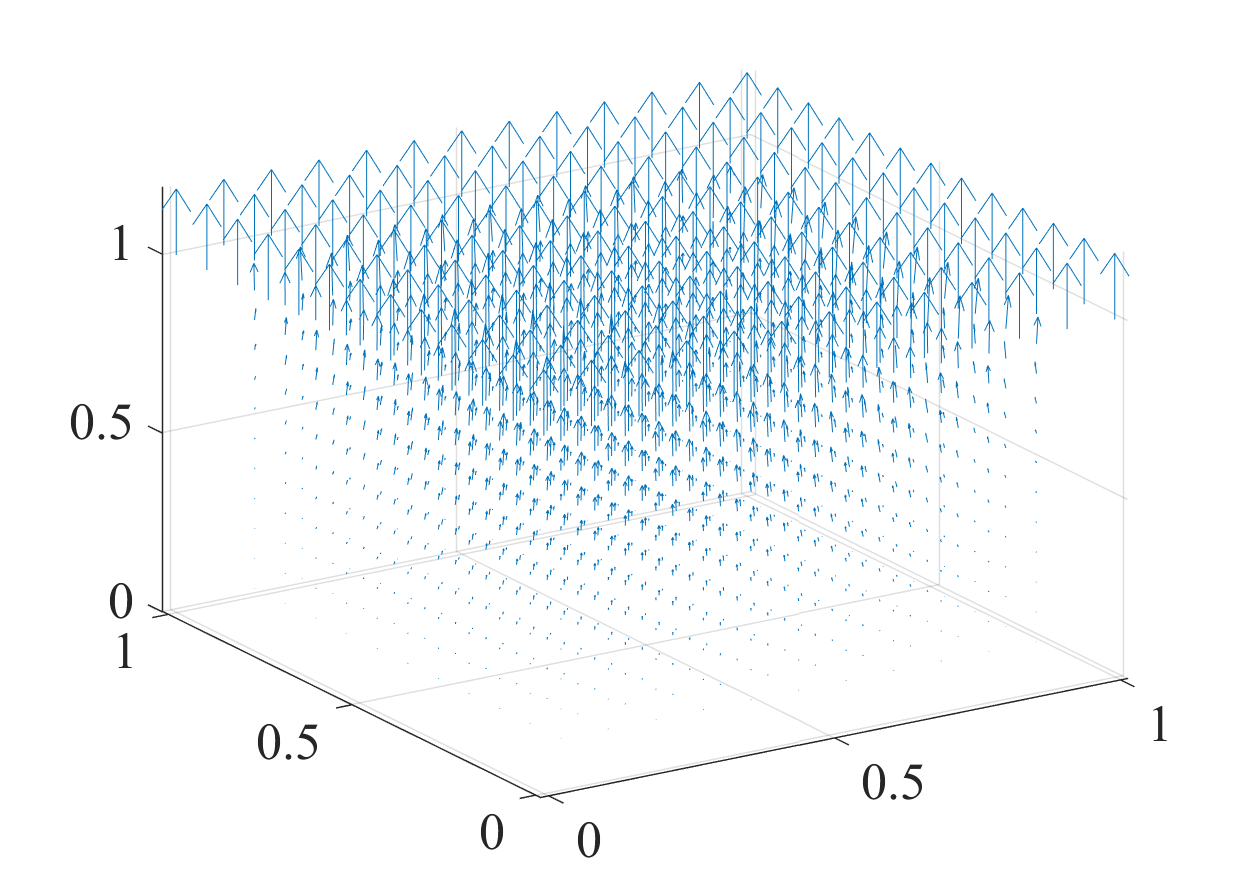}}
\caption{\label{fig:0}\small  The micro-rotation field (a) and (b) of flow by the Scheme 5.1 and Scheme 5.3.}
\end{figure}

%%%%%%%%%%%%%%%%%%%%%%%%%%%%%%%%%%%%%%%%%%%%%%%%%%%%%%%%%%%%%%%%%%%%%%%%%%%%%%%%%%%%%%%%%%%%%%%%%%%%%%%%%%%%%%%%%%%%%%

%%%%%%%%%%%%%%%%%%%%%%%%%%%%%%%%%%%%%%%%%%%%%%%%%%%%%%%%%%%%%%%%%%%%%%%%%%%%%%%%%%%%%%%%%%%%%%%%%%%%%%%%%%%%%%%%%%%%%
%
%
%
%%%%%%%%%%%%%%%%%%%%%%%%%%%%%%%%%%%%%%%%%%%%%%%%%%%%%%%%%%%%%%%%%%%%%%%%%%%%%%%%%%%%%%%%%%%%%%%%%%%%%%%%%%%%%%%%%%%%%%

\section*{\label{Sec6}Conclusions}

In this paper, we have studied some unconditionally energy stable
 numerical schemes for solving the nonstationary incompressible $3$D magneto-micropolar problem.
 The first scheme is comprised of the backward Euler semi-implicit discretization
in time and conforming finite element/stabilized finite
element in space.
The second one is based on Crank-Nicolson discretization
in time and extrapolated treatment of the nonlinear terms such that skew-symmetry
properties are retained.
 We prove that the proposed schemes are
unconditionally stable. Some optimal error estimates are derived. Moreover, we establish some
fully discrete first-order decoupled numerical schemes.
Numerical tests are given to illustrate the accuracy, efficiency, and long-time stability
of our proposed schemes.

\section*{\label{Sec6}Declarations}

\textbf{Conflict of interest}
The authors declare no competing interests.

%\begin{acknowledgements}
%If you'd like to thank anyone, place your comments here
%and remove the percent signs.
%\end{acknowledgements}

% Authors must disclose all relationships or interests that
% could have direct or potential influence or impart bias on
% the work:
%
% \section*{Conflict of interest}
%
% The authors declare that they have no conflict of interest.

% BibTeX users please use one of
%\bibliographystyle{spbasic}      % basic style, author-year citations
%\bibliographystyle{spmpsci}      % mathematics and physical sciences
%\bibliographystyle{spphys}       % APS-like style for physics
%\bibliography{}   % name your BibTeX data base

% Non-BibTeX users please use

\end{document}